# CONSISTENCY OF BAYES ESTIMATORS OF A BINARY REGRESSION FUNCTION[1]

By Marc Coram and Steven P. Lalley

*University of Chicago*

When do nonparametric Bayesian procedures "overfit"? To shed light on this question, we consider a binary regression problem in detail and establish frequentist consistency for a certain class of Bayes procedures based on hierarchical priors, called *uniform mixture priors*. These are defined as follows: let $\nu$ be any probability distribution on the nonnegative integers. To sample a function $f$ from the prior $\pi^\nu$, first sample $m$ from $\nu$ and then sample $f$ uniformly from the set of step functions from $[0, 1]$ into $[0, 1]$ that have exactly $m$ jumps (i.e., sample all $m$ jump locations and $m + 1$ function values independently and uniformly). The main result states that if a data-stream is generated according to any fixed, measurable binary-regression function $f_0 \not\equiv 1/2$, then frequentist consistency obtains: that is, for any $\nu$ with infinite support, the posterior of $\pi^\nu$ concentrates on any $L^1$ neighborhood of $f_0$. Solution of an associated large-deviations problem is central to the consistency proof.

## 1. Introduction.

### 1.1. *Consistency of Bayes procedures.*

It has been known since the work of Freedman [7] that Bayesian procedures may fail to be consistent in the frequentist sense: For estimating a probability density on the natural numbers, Freedman exhibited a prior that assigns positive mass to every open set of possible densities, but for which the posterior is consistent only at a set of the first category. Freedman's example is neither pathological nor rare: for other instances, see [4, 8, 10] and the references therein.

*Frequentist consistency* of a Bayes procedure here will mean that the posterior probability of each neighborhood of the true parameter tends to 1.

Received December 2004; revised September 2005.

[1]Supported in part by NSF Grant DMS-04-05102.

AMS 2000 subject classifications. 62A15, 62E20.

*Key words and phrases.* Consistency, Bayes procedure, binary regression, large deviations, subadditivity.

The choice of topology may be critical: For consistency in the weak topology on measures, it is generally enough that the prior should place positive mass on every Kullback–Leibler neighborhood of the true parameter [22], but for consistency in stronger topologies, more stringent requirements on the prior are needed—see, for example, [1, 9, 23]. Roughly, these demand not only that the prior charge Kullback–Leibler neighborhoods of the true parameter, but also that it not be overly diffuse, as this can lead to overfitting. Unfortunately, it appears that in certain nonparametric function estimation problems, the general formulation of this latter requirement for consistency in [1] is far too stringent (see the discussion in Section 1.3 below), as it rules out large classes of useful priors for which the corresponding Bayes procedures are in fact consistent.

1.2. *Binary regression.* The purpose of this paper is to examine in detail the consistency properties of Bayes procedures based on certain hierarchical priors in a nonparametric regression problem. For mathematical simplicity, we shall work in the setting of *binary regression*, with covariates valued in the unit interval $[0, 1]$, and we shall limit consideration to the uniform mixture priors defined below. The approach we develop can, however, be adapted to a variety of function estimation problems in one dimension (and perhaps in higher dimensions as well) and to other classes of hierarchical priors. In Section 1.4 below we provide a brief template of the approach.

Consistency of Bayes procedures in binary regression has been studied previously by Diaconis and Freedman [5, 6] for a class of priors—suggested by de Finetti—that are supported by the set of step functions with discontinuities at dyadic rationals. The use of such priors may be quite reasonable in circumstances where the covariate is actually an encoding (via binary expansion) of an infinite sequence of binary covariates. However, in applications where the numerical value of the covariate represents a real physical variable, the restriction to step functions with discontinuities only at dyadics is highly unnatural; and simulations show that when the regression function is continuous, the concentration of the posterior may be quite slow.

Coram [3] proposed another class of priors, which we shall call *uniform mixture priors*, on step functions. These are at once mathematically natural, allow computationally efficient simulation of posteriors, and appear to have much more favorable concentration properties for data generated by continuous binary regression functions than do the Diaconis–Freedman priors. In simulation experiments [3] the posterior mean of a uniform mixture prior had noticeably smaller MSE on average than CART estimates. Its performance was similar to bagged CART, but with slightly smaller MSE on average. See Figure 1 for an example.

The uniform mixture priors $\pi^\nu$, like those of Diaconis and Freedman, are hierarchical priors parametrized by probability distributions $\nu$ on the nonnegative integers. A random step function with distribution $\pi^{\nu}$ can be obtained as follows: (1) Choose a random integer $M$ with distribution $\nu$. (2) Given that $M = m$, choose $m$ points $u_i$ at random in $[0, 1]$ according to the uniform distribution: these are the points of discontinuity of the step function. (3) Given $M = m$ and the discontinuities $u_i$, choose the $m + 1$ step heights $w_j$ by sampling again from the uniform distribution. The uniform sampling in steps (2)–(3) allows for easy and highly efficient Metropolis–Hastings simulations of posteriors; the uniform distribution could be replaced by other distributions in either step, at the expense of some efficiency in posterior simulations (and our main theoretical results could easily be extended to such priors), but we see no compelling reason to discuss such generalizations in detail.

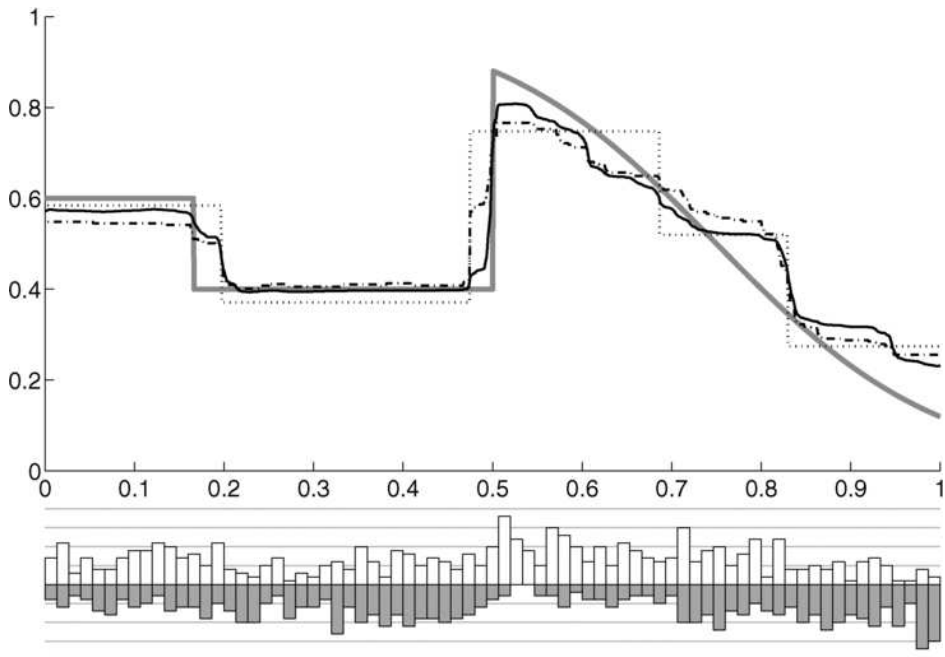

Fig. 1. *Simulation example: The data is simulated with 1024 random x-values and Bernoulli y's whose success probability is the true binary regression function, the thick gray curve. The thin solid curve is the posterior mean of the uniform mixture prior with $\nu$ chosen to be Geometric($\frac{1}{2}$). For comparison, the dotted line is cross-validated CART and the dash-dotted line is bagged CART. The white and gray histograms at the bottom show the raw data.*

Let $f$ be a binary regression function on $[0, 1]$, that is, a Borel-measurable function $f : [0, 1] \to [0, 1]$. We shall assume that under $P_f$ the data $(X_i, Y_i)$ are i.i.d. random vectors, with $X_i$ uniformly distributed on $[0, 1]$ and $Y_i$, given $X_i = x$, is Bernoulli-$f(x)$. (Our main result would also hold if the covariate distribution were not uniform but any other distribution giving positive mass to all intervals of positive length.) Let $Q^\nu = \int P_f \, d\pi^\nu$, and denote by $Q^\nu(\cdot | \mathcal{F}_n)$ the posterior distribution on step functions given the first $n$ observations $(X_i, Y_i)$ [more precisely, the regular version of the conditional distribution defined by (17) below]. The main result of the paper is as follows.

THEOREM 1. *Assume that the hierarchy prior $\nu$ is not supported by a finite subset of the integers. Then for every binary regression function $f \not\equiv \frac{1}{2}$, the $Q^\nu$-Bayes procedure is $L^1$-consistent at $f$, that is, for every $\varepsilon > 0$,*

$$(1) \qquad \lim_{n \to \infty} P_f \{ Q^\nu(\{g : \|g - f\|_1 > \varepsilon\} | \mathcal{F}_n) > \varepsilon \} = 0.$$

The restriction $f \not\equiv 1/2$ arises for precisely the same reason as in [6], namely, that this exceptional function is the prior mean of the regression function. See [6] for further discussion.

Theorem 1 implies that the uniform mixture priors enjoy the same consistency as do the Diaconis–Freedman priors [6]. This is not exactly unexpected, but neither should it be considered a priori obvious—as the proof will show, there are substantial differences between the uniform mixture priors and those of Diaconis and Freedman: In particular, since the uniform mixture priors allow the step-function discontinuities to arrange themselves in favorable (but atypical for uniform samples) configurations *vis-a-vis* the data, the danger of overfitting would seem, at least a priori, greater than for the Diaconis–Freedman priors. The bulk of the proof (Sections 4–5) will be devoted to showing that such overfitting does not occur, except possibly when $f \equiv 1/2$.

Theorem 1 asserts only weak convergence (i.e., in $P_f$-probability). In fact, the arguments can be extended to establish almost sure convergence, but at the expense of added complication. This would involve replacing the subadditive WLLN in Appendix A by a corresponding strong law, and modifying those arguments in Section 4 used to verify the hypotheses of the subadditive WLLN.

1.3. *Relation to other work.* There is a substantial literature on the consistency of Bayes procedures, much of it devoted to establishing sufficient conditions. See, for a start, [1, 10, 22, 25, 26] and the references therein. Certain of the sufficient conditions developed in these papers apply, at least in principle, to hierarchical priors of the type considered here. Unfortunately, these conditions require that the prior be highly concentrated on low-complexity models. For instance, the main result of [1] would require for the uniform mixture priors that the hierarchy prior $\nu$ satisfy

$$\sum_{k \geq m} \nu_k \leq m^{-mC}$$

for some $C > 0$ (see [3]). Recent results of Walker [25] improve those of [1], but, for a uniform mixture, priors evidently still require that $\nu$ have an exponentially decaying tail ([25], Section 6.3). Such restrictions on the tail of the hierarchy prior certainly prevent the accumulation of posterior mass on models that are overfit, but at the possible cost of having the posterior favor models that are *underfit*. Preliminary analysis seems to indicate that, when the true regression function is smooth, more rapid posterior concentration takes place when the hierarchy prior has a rather long tail. One objective of this paper is to show that, in at least one model of real statistical interest, the problem of finding the *right* conditions for frequentist consistency of Bayesian procedures requires a careful analysis of an associated large deviations problem.

1.4. *Overfitting and large deviations problems.* The possibility of overfitting is tied up with a certain *large deviations* problem connected to the model: this is the most interesting mathematical feature of the paper. (A similar large deviations problem occurs in [6], but there it reduces easily to the classical Cramér LD theorem for sums of i.i.d. random variables.) Roughly, we will show in Section 4 that, as the sample size $n \to \infty$, the posterior probability of the set of step functions with more than $\alpha n$ discontinuities decays like $e^{n\psi(\alpha)}$, where $\psi(\alpha) < 0$. Then, in Section 5, we will show that $\psi(\alpha)$ is uniquely maximized at $\alpha \to 0$; this will imply that, for large sample size $n$, most of the posterior mass is concentrated on step functions with a small number of discontinuities relative to $n$. Concentration of the posterior in $L^1$-neighborhoods of the true regression function will then follow by routine arguments—see Section 3.

We expect (and hope to show in a subsequent paper) that in a variety of problems, for certain classes of hierarchical priors, the critical determinant of the consistency of Bayes procedures will prove to be the rate functions in associated large deviations problems. The template of the analysis is as follows: Let

$$
(2) \qquad \pi = \pi^\nu = \sum_{m=0}^{\infty} \nu_m \pi_m
$$

be a hierarchical prior obtained by mixing priors $\pi_m$ of "complexity" $m$. Let $Q$ and $Q_m$ be the probability distributions on the space of data sequences gotten by mixing with respect to $\pi$ and $\pi_m$, respectively; and let $Q(\cdot | \mathcal{F}_n)$ and $Q_m(\cdot | \mathcal{F}_n)$ be the corresponding posterior distributions given the information in the $\sigma$-field $\mathcal{F}_n$. Then by Bayes' formula,

$$
(3) \qquad Q(\cdot | \mathcal{F}_n) = \left\{ \sum_{m=0}^{\infty} \nu_m Z_{m,n} Q_m(\cdot | \mathcal{F}_n) \right\} \Big/ \left\{ \sum_{m=0}^{\infty} \nu_m Z_{m,n} \right\},
$$

where $Z_{m,n}$ are the *predictive probabilities* for the data in $\mathcal{F}_n$ based on the model $Q_m$ (see Section 2.2 for more detail in the binary regression problem). This formula makes apparent that the relative sizes of the predictive probabilities $Z_{m,n}$ determine where the mass in the posterior $Q^\nu(\cdot|\mathcal{F}_n)$ is concentrated. The *large deviations problem* is to show that as $m, n \to \infty$ in such a way that $m/n \to \alpha$,

$$(4) \qquad Z_{m,n}^{1/n} \longrightarrow \exp\{\psi(\alpha)\}$$

in $P_f$-probability, for an appropriate nonrandom *rate function* $\psi(\alpha)$, and to show that $\psi(\alpha)$ is uniquely maximized at $\alpha = 0$. This, when true, will imply that most of the posterior mass will be concentrated on models with small complexity $m$ relative to the sample size $n$, where overfitting does *not* occur.

1.5. *Choice of topology.* The use of the $L^1$-metric (equivalently, any $L^p$-metric, $0 \le p < \infty$) in measuring posterior concentration, as in (1), although in many ways natural, may not always be appropriate. Posterior concentration relative to the $L^1$-metric justifies confidence that, for a new random sample of individuals with covariates uniformly distributed on $[0, 1]$, the responses will be reasonably well-predicted by regression function samples from the posterior, but it would not justify similar confidence for a random sample of individuals all with covariate (say) $x = 0.47$. For this, posterior concentration in the sup-norm metric would be required. We do not yet know if consistency holds in the sup-norm metric, for either the uniform mixture priors or the Diaconis–Freedman priors, even for smooth $f$; but we conjecture that it does.

## 2. Preliminaries.

2.1. *Data.* A (binary) *regression function* is a Borel measurable function $f : J \to [0, 1]$, where $J$ is an interval. Most often the interval $J$ will be the unit interval. For each binary regression function $f$, let $P_f$ be a probability measure on a measurable space supporting a *data stream* $\{(X_n, Y_n)\}_{n \ge 1}$ such that under $P_f$ the random variables $X_n$ are i.i.d. Uniform-$[0, 1]$ and, conditional on $\sigma(\{X_n\}_{n \ge 1})$, the random variables $Y_n$ are independent Bernoullis with conditional means

$$(5) \qquad E_f(Y_n|\sigma(\{X_m\}_{m \ge 1})) = f(X_n).$$

(In several arguments below it will be necessary to consider alternative distributions $F$ for the covariates $X_n$. In such cases we shall adopt the convention of adding the subscript $F$ to relevant quantities; thus, for instance, $P_{f,F}$ will denote a probability distribution under which the covariates $X_n$ are i.i.d. $F$, and the conditional distribution of the responses $Y_n$ is the same as under $P_f$.) We shall assume when necessary that probability spaces support additional independent streams of uniform and exponential r.v.s (and thus also Poisson processes), so that auxiliary randomization is possible. Generic data sets [values of the first $n$ pairs $(x_i, y_i)$] will be denoted $(\mathbf{x}, \mathbf{y})$ or $(\mathbf{x}, \mathbf{y})_n$ to emphasize the sample size; the corresponding random vectors will be denoted by the matching upper case letters $(\mathbf{X}, \mathbf{Y})$. For any data set $(\mathbf{x}, \mathbf{y})$ and any interval $J \subset [0, 1]$, the number of successes $(y_i = 1)$, failures $(y_i = 0)$ and the total number of data points with covariate $x_i \in J$ will be denoted by

$$N^S(J), N^F(J) \quad \text{and} \quad N(J) = N^S(J) + N^F(J).$$

In certain comparison arguments, it will be convenient to have data streams for different regression functions defined on a common probability space $(\Omega, \mathcal{F}, P)$. This may be accomplished by the usual device: Let $\{X_n\}_{n \geq 1}$ and $\{V_n\}_{n \geq 1}$ be independent, identically distributed Uniform-$[0, 1]$ random variables, and set

$$(6) \qquad Y_n^f = \mathbb{1}\{V_n \leq f(X_n)\}.$$

2.2. *Priors on regression functions.* The prior distributions on regression functions considered in this paper are probability measures on the set of step functions with finitely many discontinuities. Points of discontinuity, or *split points*, of step functions will be denoted by $u_i$, and step heights by $w_i$. Each vector $\mathbf{u} = (u_1, u_2, \ldots, u_m)$ of split points induces a partition of the unit interval into $m + 1$ subintervals (or *cells*) $J_i = J_i(\mathbf{u})$. [Note: We do not assume that split point vectors $(u_1, u_2, \ldots, u_m)$ are ordered.] Denote by $\pi_{\mathbf{u}}$ the probability measure on step functions with discontinuities $u_i$ that makes the step height random variables $W_i$ [i.e., the values $w_i$ on the intervals $J_i(\mathbf{u})$] independent and uniformly distributed on $[0, 1]$. For each nonnegative integer $m$, define $\pi_m$ to be the uniform mixture of the measures $\pi_{\mathbf{u}}$ over all split point vectors $\mathbf{u}$ of length $m$, that is,

$$(7) \qquad \pi_m = \int_{\mathbf{u} \in (0,1)^m} \pi_{\mathbf{u}} \, d\mathbf{u}.$$

It is, of course, possible to mix against distributions $G$ other than the uniform, and in some arguments it will be necessary for us to do so: in such cases (see, e.g., Section 2.3) we shall use additional subscripts $G$ on various objects to indicate that the split point vectors are gotten by sampling from $G$ instead of from the uniform distribution $U$. The priors of primary interest—those considered in Theorem 1 and equation (2)—are mixtures $\pi^\nu = \sum_m \nu_m \pi_m$ of the measures $\pi_m$ against *hierarchy priors* $\nu$ on the nonnegative integers $\mathbb{N}$. Unless otherwise stated, assume throughout that the hierarchy prior $\nu$ is not supported by a finite subset of $\mathbb{N}$.

Each of the probability measures $\pi_{\mathbf{u}}$, $\pi_m$ and $\pi^\nu$ induces a corresponding probability measure on the space of data sequences by mixing:

$$\text{(8)} \qquad\qquad Q_{\mathbf{u}} = \int P_f \, d\pi_{\mathbf{u}}(f),$$

$$\text{(9)} \qquad\qquad Q_m = \int P_f \, d\pi_m(f)$$

and

$$\text{(10)} \qquad\qquad Q^\nu = \int P_f \, d\pi^\nu(f).$$

Observe that $Q_m$ is the uniform mixture of the measures $Q_{\mathbf{u}}$ over split point vectors $\mathbf{u}$ of size $m$, and $Q^\nu$ is the $\nu$-mixture of the measures $Q_m$.

For any data sample $(\mathbf{x}, \mathbf{y})$, the *posterior distribution* $Q(\cdot|(\mathbf{x}, \mathbf{y}))$ under any of the measures $Q_{\mathbf{u}}$, $Q_m$ or $Q^\nu$ is the conditional distribution on the set of step functions given that $(\mathbf{X}, \mathbf{Y}) = (\mathbf{x}, \mathbf{y})$. The posterior distribution $Q_{\mathbf{u}}(\cdot|(\mathbf{x}, \mathbf{y}))$ can be explicitly calculated: it is the distribution that makes the step height r.v.s $W_i$ independent, with Beta-$(N_i^S, N_i^F)$ distributions, where $N_i^S = N^S(J_i(\mathbf{u}))$ and $N_i^F = N^F(J_i(\mathbf{u}))$ are the success/failure counts in the intervals $J_i$ of the partition induced by $\mathbf{u}$. Thus, the joint density of the step heights (relative to product Lebesgue measure on the cube $[0,1]^{m+1}$ is

$$\text{(11)} \qquad q_{\mathbf{u}}(\mathbf{w}|(\mathbf{x}, \mathbf{y})) = Z_{\mathbf{u}}(\mathbf{x}, \mathbf{y})^{-1} \prod_{i=0}^m w_i^{N_i^S} (1 - w_i)^{N_i^F},$$

where the normalizing constant $Z_{\mathbf{u}}(\mathbf{x}, \mathbf{y})$, henceforth called the $Q_{\mathbf{u}}$-*predictive probability* for the data sample $(\mathbf{x}, \mathbf{y})$, is given by

$$\text{(12)} \quad Z_{\mathbf{u}}(\mathbf{x}, \mathbf{y}) = \int_{\mathbf{w} \in [0,1]^{m+1}} \prod_{i=0}^m w_i^{N_i^S} (1 - w_i)^{N_i^F} \, d\mathbf{w} = \prod_{i=0}^m B(N_i^S, N_i^F)$$

and

$$\text{(13)} \qquad\qquad B(m, n) = \left\{ (m + n + 1) \binom{m + n}{m} \right\}^{-1}.$$

(This is not the usual convention for the arguments of the beta function, but will save us from a needless proliferation of +1s.) The posterior distributions $Q_m(\cdot|(\mathbf{x}, \mathbf{y}))$ and corresponding predictive probabilities $Z_m(\mathbf{x}, \mathbf{y})$ are related to $Q_{\mathbf{u}}(\cdot|(\mathbf{x}, \mathbf{y}))$ and $Z_{\mathbf{u}}(\mathbf{x}, \mathbf{y})$ as follows:

$$\text{(14)} \quad Q_m(\cdot|(\mathbf{x}, \mathbf{y})) = \left\{ \int_{\mathbf{u} \in [0,1]^m} Q_{\mathbf{u}}(\cdot|(\mathbf{x}, \mathbf{y})) Z_{\mathbf{u}}(\mathbf{x}, \mathbf{y}) \, d(\mathbf{u}) \right\} \Big/ Z_m(\mathbf{x}, \mathbf{y}),$$

where

$$
\begin{aligned}
Z_m(\mathbf{x}, \mathbf{y}) &= \int_{\mathbf{u} \in (0,1)^m} Z_{\mathbf{u}}(\mathbf{x}, \mathbf{y}) \, d(\mathbf{u}) \\
&= \int_{\mathbf{u} \in (0,1)^m} \prod_{i=0}^{m} B(N_i^S, N_i^F) \, d(\mathbf{u}).
\end{aligned}
$$

(15)

(Note: The dependence of the integrand on $\mathbf{u}$, via the values of the success/failure counts $N_i^S, N_i^F$, is suppressed.) In general, the last integral cannot be evaluated in closed form, unlike the integral (12) that defines the $Q_{\mathbf{u}}$-predictive probabilities. This, as we shall see in Sections 4–5, will make the mathematical analysis of the posteriors considerably more difficult than the corresponding analysis for Diaconis–Freedman priors.

Note for future reference (Section 5) that the predictive probabilities $Z_m$ are related to likelihood ratios $dQ_m/dP_f$: In particular, when $f \equiv p$ is constant,

$$
Z_m((\mathbf{X}, \mathbf{Y})_n) = p^{N^S}(1-p)^{N^F} \left( \frac{dQ_m}{dP_p} \right)_{\mathcal{F}_n},
$$

(16)

where $\mathcal{F}_N$ is the $\sigma$-algebra generated by the first $n$ data points, and $N^S, N^F$ are the numbers of successes and failures in the entire data set $(\mathbf{X}, \mathbf{Y})_n$. Finally, observe that the posterior distribution $Q^\nu(\cdot | (\mathbf{x}, \mathbf{y}))$ is related to the posterior distributions $Q_m(\cdot | (\mathbf{x}, \mathbf{y}))$ by Bayes' formula,

(17) $\quad Q^\nu(\cdot | (\mathbf{x}, \mathbf{y})) = \left\{ \sum_{m=0}^{\infty} \nu_m Z_m(\mathbf{x}, \mathbf{y}) Q_m(\cdot | (\mathbf{x}, \mathbf{y})) \right\} \Big/ \left\{ \sum_{m=0}^{\infty} \nu_m Z_m(\mathbf{x}, \mathbf{y}) \right\}.$

The goal of Sections 4–5 will be to show that, for large samples $(\mathbf{X}, \mathbf{Y})_n$, under $P_f$ the predictive probabilities $Z_{\alpha n}((\mathbf{X}, \mathbf{Y})_n)$ are of smaller exponential magnitude for $\alpha > 0$ than for $\alpha = 0$. This will imply that the posterior concentrates in the region $m \ll n$, where the number of split points is small compared to the number of data points.

CAUTION. Note that $\pi_m$ and $\pi_{\mathbf{u}}$ have different meanings, as do $Z_m$ and $Z_{\mathbf{u}}$, and $Q_{\mathbf{u}}$ and $Q_m$. The reader should have no difficulty discerning the proper meaning by context or careful examination of the fonts.

2.3. *Transformations of the covariates.* We have assumed that under $P_f$ the covariates $X_n$ are uniformly distributed on $[0,1]$; and, in constructing the priors $\pi_m$ and $\pi^\nu$, we have used uniform mixtures on the locations of the split points $u_i$. This is only for convenience: the covariate space could be relabeled by any homeomorphism without changing the nature of the estimation problem. Thus, if the data sample $(\mathbf{x}, \mathbf{y})$ were changed to $(G^{-1}\mathbf{x}, \mathbf{y})$, where $G$ is a continuous, strictly increasing c.d.f., and if $G$-mixtures rather than uniform mixtures were used in building the priors, then the predictive probabilities would be unchanged:

$$(18) \qquad Z_{m,G}(G^{-1}\mathbf{x}, \mathbf{y}) = Z_m(\mathbf{x}, \mathbf{y}),$$

where $Z_m = Z_{m,G}$ are the predictive probabilities for the transformed data relative to priors $\pi_G^\nu$ built using $G$-mixtures instead of uniform mixtures. (This follows directly from the transformation formula for integrals.)

In a number of arguments below it will be necessary to consider data streams $(\mathbf{X}, \mathbf{Y})$ distributed according to $P_{f,F}$, where $F$ is a distribution other than the uniform distribution on $[0, 1]$. First, if $(\mathbf{x}, \mathbf{y})_n$ is a data sample of size $n$, with covariates $x_i \in [0, 1]$, and if $G$ is the c.d.f. of the uniform distribution on the interval $[0, n]$ (so that $G^{-1}$ is just multiplication by $n$), then applying the transformation $G^{-1}$ has the effect of standardizing the spacings between data points and between split points. This will be used in Section 5. Second, if the data stream $\{(X_n, Y_n)\}_{n \geq 1}$ is subjected to *thinning*—for instance, remove a data point $(X_n, Y_n)$ from the stream with probability $\rho(X_n)$ depending on the value of the covariate $X_n$—then the resulting thinned data stream will be distributed according to $P_{f,F}$, where $F$ has density proportional to $1 - \rho$. The comparison arguments in Section 4.6 below will rely on this fact.

2.4. *Self-similarity.* The key to analyzing the predictive probability (15) is that the integral in (15) has an approximately self-similar structure: it almost (but not exactly!) factors as the product of two integrals each of the same form, one over the data and split points in $[0, 1/2]$, the other over $(1/2, 1]$:

$$(19) \qquad Z_m(\mathbf{X}, \mathbf{Y}) \approx Z_{m/2,G_0}(\mathbf{X}', \mathbf{Y}') Z_{m/2,G_1}(\mathbf{X}'', \mathbf{Y}'').$$

Here $G_0, G_1$ are the uniform distributions on $[0, 0.5]$ and $[0.5, 1]$, respectively, and $(\mathbf{X}', \mathbf{Y}')$ and $(\mathbf{X}'', \mathbf{Y}'')$ are the subsets of the data set $(\mathbf{X}, \mathbf{Y})$ with covariates in $[0, 0.5]$ and $[0.5, 1]$, respectively. Unfortunately, the factorization is not exact, for two reasons: (i) the split point vectors $\mathbf{u}$ in the integral (15) do not necessarily include a split point at $1/2$; and (ii) the number of split points in $(0, 1/2)$ is not exactly $m/2$. Nevertheless, when $m$ is large, most split point vectors $\mathbf{u}$ will include a split very near $1/2$, and will put about $m/2$ splits in each of $(0, 1/2)$ and $(1/2, 1)$, and so it is not unreasonable to expect that (19) should hold approximately.

Consider the two factors $Z', Z''$ in the approximation (19). If the sample size $n$ is large, then each of the subsamples $(\mathbf{X}', \mathbf{Y}')$ and $(\mathbf{X}'', \mathbf{Y}'')$ should contain about $n/2$ points. Furthermore, if the true regression function $f \equiv p$ is constant, then, under $P_p$ each of the subsamples should, after covariate transformation by $x \mapsto 2x \bmod 1$, be distributed as a sample of size (about) $n/2$ under $P_p$. Therefore, by (18), the factors in (19) are under $P_p$ independent random variables, each with the same distribution as $Z_{m/2}((\mathbf{X}, \mathbf{Y})_{n/2})$ under $P_p$.

Iteration of this factorization exhibits the predictive probability $Z_m(\mathbf{X}, \mathbf{Y})$ approximately as a product of a large number of independent, identically distributed factors. (Note, though, that the errors in these approximations may accumulate exponentially in the number of iterations; this will be handled by the use of subadditivity arguments in Section 4.) In essence, the exponential decay (4) follows from this approximate product representation.

2.5. *Beta function and Beta distributions.* Because the posterior distributions (11) of the step height random variables are Beta distributions, certain elementary properties of these distributions and the corresponding normalizing constants $B(n, m)$ will play an important role in the analysis. The behavior of the Beta function for large arguments is well understood, and easily deduced from Stirling's formula. Similarly, the asymptotic behavior of the Beta distributions follows from the fact that these are the distributions of uniform order statistics:

BETA CONCENTRATION PROPERTY. *For each $\varepsilon > 0$, there exists $k(\varepsilon) < \infty$ such that, for all index pairs $(m, n)$ with weight $m + n > k(\varepsilon)$, (a) the Beta-$(m, n)$ distribution puts all but at most $\varepsilon$ of its mass within $\varepsilon$ of $m/(m + n)$; and* (b) *the normalization constant $B(m, n)$ satisfies*

$$(20) \qquad \left| \frac{\log B(m, n)}{m + n} + H\left(\frac{m}{m + n}\right) \right| < \varepsilon,$$

*where $H(x)$ is the Shannon entropy, defined by*

$$(21) \qquad H(x) = -x \log x - (1 - x) \log(1 - x).$$

Note that the binomial coefficient in (13) is bounded above by $2^{m+n}$, so it follows that $B(m, n) \geq 4^{-m-n}$. Thus, by (15), for any data sample $(\mathbf{x}, \mathbf{y})$ of size $n$,

$$(22) \qquad Z_m(\mathbf{x}, \mathbf{y}) \geq 4^{-n}.$$

Some of the arguments in Section 4 will require an estimate of the effect on the integral (15) of adding another split point. This breaks one of the intervals $J_i$ into two, leaving all of the others unchanged, and so the effect on the integrand in (15) is that one of the factors $B(N_i^S, N_i^F)$ is replaced by a product of two factors $B(N_L^S, N_L^F)B(N_R^S, N_R^F)$, where the cell counts satisfy

$$N_L^S + N_R^S = N_i^S$$

and

$$N_L^F + N_R^F = N_i^F.$$

The following inequality, which is easily deduced from (13), shows that the multiplicative error made in this replacement is bounded by the overall sample size:

$$(23) \qquad \frac{B(N_i^S, N_i^F)}{B(N_L^S, N_L^F) B(N_R^S, N_R^F)} \leq N_i^S + N_i^F.$$

2.6. *The entropy functional.* We will show, in Section 3 below, that in the "Middle Zone," where the number of split points is large but small compared to the number $n$ of data points, the predictive probability decays at a precise exponential rate as $n \to \infty$. The rate is the negative of the *entropy functional* $H(f)$, defined by

$$(24) \qquad H(f) = \int_0^1 H(f(x)) \, dx,$$

where $H(x)$ for $x \in (0, 1)$ is the Shannon entropy defined by (21) above. The Shannon entropy function $H(x)$ is uniformly continuous and strictly concave on $[0, 1]$, with second derivative bounded away from 0; it is strictly positive except at the endpoints, 0 and 1; and it attains a maximum value of $\log 2$ at $x = 1/2$. The entropy functional $H(f)$ enjoys similar properties:

ENTROPY CONTINUITY PROPERTY. *For each $\varepsilon > 0$, there exists $\delta > 0$ so that*

$$(25) \qquad \|f - g\|_1 < \delta \quad \implies \quad |H(f) - H(g)| < \varepsilon.$$

ENTROPY CONCAVITY PROPERTY. *Let $f$ and $g$ be binary regression functions such that $g$ is an averaged version of $f$ in the following sense: There exist finitely many pairwise disjoint Borel sets $B_i$ such that $\{x : g(x) \neq f(x)\} = \bigcup_i B_i$, and for each $i$ such that $|B_i| > 0$,*

$$(26) \qquad g(x) = \int_{B_i} f(y) \, dy / |B_i| \qquad \forall x \in B_i.$$

*Then*

$$(27) \qquad H(g) - H(f) \geq -\left(\max_{0 < p < 1} H''(p)\right) \|f - g\|_1 / 2.$$

*Hence, $H(g) \geq H(f)$ with strict inequality unless $f = g$ a.e.*

PROOF. The Continuity property (25) follows from the uniform continuity of the Shannon function $H(p)$ by an elementary argument. The Concavity property follows from the Jensen inequality, and the "uniform" strengthening (27) from the fact that $H''(p)$ is bounded away from zero on $[0, 1]$. To see this, let $B = \{x : f(x) \neq g(x)\}$, and let $-C < 0$ be an upper bound for $H''(p)$. The hypothesis (26) implies that $g$ is constant on $B_i$, and equal to the average $a_i$ of $f$ on this set. By Taylor's theorem, on $B_i$,

$$H(f(x)) - H(a_i) = H'(a_i)(f(x) - a_i) + H''(\zeta(x))(f(x) - a_i)^2/2,$$

where $\zeta(x)$ is intermediate between $a_i$ and $f(x)$. Hence,

$$\begin{aligned} H(f) - H(g) &= \int_B H''(\zeta(x))(f(x) - g(x))^2 \, dx/2 \\ &\leq -C\|f - g\|_2^2/2 \\ &\leq -C\|f - g\|_1/2, \end{aligned}$$

the last inequality following because $0 \leq f, g \leq 1$. $\quad\square$

The Continuity and Concavity properties will be used principally to estimate the entropies of step functions $g$ near $f$. In particular, they allow entropy comparisons between a binary regression function $f$ and step functions obtained by averaging $f$ on the intervals of a partition. Let $\mathbf{u}$ be a vector of split points, and let $J_i = J_i(\mathbf{u})$ be the intervals in the partition of $[0, 1]$ induced by $\mathbf{u}$. For each binary regression function $f$, define $\bar{f}_{\mathbf{u}}$ to be the step function whose value on each interval $J_i(\mathbf{u})$ is the mean value $\int_{J_i} f/|J_i|$ of $f$ on that interval. Then by the Concavity property,

$$(28) \qquad H(\bar{f}_{\mathbf{u}}) \geq H(f),$$

with strict inequality unless $f = \bar{f}_{\mathbf{u}}$ a.e. Moreover, the difference is small if and only if $f$ and $\bar{f}_{\mathbf{u}}$ are close in $L^1$. This will be the case if all intervals $J_i$ of the partition are small:

LEMMA 1. *For each binary regression function $f$ and each $\varepsilon > 0$, there exists $\delta > 0$ such that if $|J_i| < \delta$ for every interval $J_i$ in the partition induced by $\mathbf{u}$, then*

$$(29) \qquad \|f - \bar{f}_{\mathbf{u}}\|_1 < \varepsilon.$$

PROOF. First, observe that the assertion is elementary for *continuous* regression functions, since continuity implies uniform continuity on $[0, 1]$. Second, recall that continuous functions are dense in $L^1[0, 1]$ by Lusin's theorem; thus, for each regression function $f$ and any $\eta > 0$, there exists a continuous function $g : [0, 1] \to [0, 1]$ such that $\|f - g\|_1 < \eta$. It then follows by the elementary inequality $\left|\int h\right| \le \int |h|$ that, for any vector $\mathbf{u}$ of split points,

$$\|\bar{f}_{\mathbf{u}} - \bar{g}_{\mathbf{u}}\|_1 < \eta.$$

Finally, use $\eta = \varepsilon/3$ and choose $\delta$ so that, for the continuous function $g$ and any $\mathbf{u}$ that induces a partition whose intervals are all of length $< \delta$,

$$\|g - \bar{g}_{\mathbf{u}}\|_1 < \eta.$$

Then by the triangle inequality for $L^1$,

$$\|f - \bar{f}_{\mathbf{u}}\|_1 \le \|g - f\|_1 + \|g - \bar{g}_{\mathbf{u}}\|_1 + \|\bar{f}_{\mathbf{u}} - \bar{g}_{\mathbf{u}}\|_1$$
$$\le 3\eta = \varepsilon. \qquad \square$$

2.7. *Empirical distributions under* $P_f$. Theorem 1 will be proved by showing (A) that, for large $n$, the posterior mass of the step functions with more than $\varepsilon n$ discontinuities is of smaller order of magnitude than that of the step functions with fewer than $\varepsilon n$ discontinuities; and (B) that the posterior mass on the latter set concentrates on those step functions that are $L^1$-close to the true regression function. Step (B) will rely on a uniform version of the law of large numbers (LLN).

Following is a suitable version of the LLN. Given a data set $(\mathbf{X}, \mathbf{Y})_n$ of size $n$ and an interval $J$, say that $J$ is $\varepsilon$-*bad* (relative to the data set) if at least one of the following inequalities holds:

$$(30) \qquad |N(J) - n|J|| \ge \varepsilon n |J|,$$

$$(31) \qquad \left| N^S(J) - n \int_J f(x)\,dx \right| \ge \varepsilon n |J|$$

or

$$(32) \qquad \left| N^F(J) - n \int_J (1 - f(x))\,dx \right| \ge \varepsilon n |J|.$$

Here $|J|$ denotes the Lebesgue measure of $J$. Given $x \in [0, 1]$, say that $x$ is $(\varepsilon, \kappa)$-*bad* (relative to the data set) if there is an $\varepsilon$-bad interval $J$ of length $|J| \ge \kappa/n$ that contains $x$. Let $\mathcal{B}_n(\varepsilon, \kappa)$ be the set of $(\varepsilon, \kappa)$-bad points relative to $(\mathbf{X}, \mathbf{Y})_n$.

PROPOSITION 2. *For any* $\varepsilon > 0$, *there exist positive constants* $\kappa, \gamma, C$ *such that, for every sample size* $n \ge 1$,

$$(33) \qquad P_f\{|\mathcal{B}_n(\varepsilon, \kappa)| \ge \varepsilon\} \le Ce^{-\gamma n}.$$

The exponential estimates, in conjunction with the Borel–Cantelli lemma, yield as a consequence a uniform *strong* law of large numbers. For the proof of Theorem 1, only a weak law is needed; however, it is no more difficult to establish the exponential bounds (33). The proof of Proposition 2 is deferred to Appendix B.

**3. Beginning and middle zones.** Following [5], we designate three asymptotic "zones" where the predictive probabilities $Z_m((\mathbf{X}, \mathbf{Y})_n)$ decay at different exponential rates. These are determined by the relative sizes of $m$, the number of discontinuities of the step functions, and $n$, the sample size. The *end zone* is the set of pairs $(m, n)$ such that $m/n \geq \varepsilon$; this zone will be analyzed in Sections 4 and 5, where we shall prove that the asymptotic decay of $Z_m((\mathbf{X}, \mathbf{Y})_n)$ is faster than in the *middle zone*, where $K \leq m \leq \varepsilon n$ for a large constant $K$. The *beginning zone* is the set of pairs $(m, n)$ for which $m \leq K$ for some large $K$. A regression function cannot be arbitrarily well approximated by step functions with a bounded number of discontinuities unless it is itself a step function, and so, as we will see, the asymptotic decay of $Z_m((\mathbf{X}, \mathbf{Y})_n)$ is generally faster in the beginning zone than in the middle zone.

In this section we analyze the beginning and middle zones, using the Beta concentration property, Lemma 1 and Proposition 2. In the beginning and middle zones, the number $m$ of split points is small compared to the number $n$ of data points, and so for typical split-point vectors $\mathbf{u}$, most intervals in the partition induced by $\mathbf{u}$ will, with high probability, contain a large number of data points. Consequently, the law of large numbers applies in these intervals: together with the Beta concentration property, it ensures that the $Q_{\mathbf{u}}$-posterior is concentrated in an $L^1$-neighborhood of $\bar{f}_{\mathbf{u}}$, and that the $Q_{\mathbf{u}}$-predictive probability is roughly $\exp\{-nH(\bar{f}_{\mathbf{u}})\}$. The next proposition makes this precise.

PROPOSITION 3. *For each $\delta > 0$, there exists $\varepsilon > 0$ such that the following is true: For all sufficiently large $n$, the $P_f$-probability is at least $1 - \delta$ that, for all $m \leq \varepsilon n$ and all split-point vectors $\mathbf{u}$ of size $m$,*

$$Q_{\mathbf{u}}(\{g : \|g - \bar{f}_{\mathbf{u}}\|_1 \geq \delta\}|(\mathbf{X}, \mathbf{Y})_n) < \delta \tag{34}$$

*and*

$$|n^{-1} \log Z_{\mathbf{u}}((\mathbf{X}, \mathbf{Y})_n) + H(\bar{f}_{\mathbf{u}})| < \delta. \tag{35}$$

PROOF. Let $J_i = J_i(\mathbf{u})$ be the intervals in the partition induced by $\mathbf{u}$. Fix $\kappa = \kappa(\delta)$ as in Proposition 2. If $\varepsilon$ is sufficiently small, then for any split-point vector $\mathbf{u}$ of size $m \leq \varepsilon n$, the union of those $J_i$ of length $\leq \kappa/n$ will have Lebesgue measure $< \delta$: this follows by a trivial counting argument. Let $B_{\mathbf{u}}$ be the union of those $J_i$ that are either of length $\leq \kappa/n$ or are $\delta$-bad [in the sense of inequalities (30)–(32)]. By Proposition 2, the $P_f$-probability of the event $G^c$ that there exists a split-point vector $\mathbf{u}$ of size $m \leq \varepsilon n$ for which the Lebesgue measure of $B_{\mathbf{u}}$ exceeds $2\delta$ is less than $\varepsilon$, for all large $n$. But on the complementary event $G$, inequality (34) must hold (with possibly different values of $\delta$) by the Beta concentration property.

For the proof of (35), recall that by (12),

$$(36) \qquad n^{-1} \log Z_{\mathbf{u}}((\mathbf{X}, \mathbf{Y})_n) = n^{-1} \sum_{i=0}^{m} \log B(N_i^S, N_i^F),$$

where $B(k, l)$ is the Beta function (using our convention for the arguments). By the Stirling approximation (20), each term of the sum for which $N_i$ is large is well approximated by $-N_i H(N_i^S / N_i)$; and for each index $i$ such that $J_i \not\subset B_{\mathbf{u}}$, this in turn is well approximated by $n|J_i| H(\bar{f}_{\mathbf{u}}(J_i))$, where $\bar{f}_{\mathbf{u}}(J_i)$ is the average of $f$ on the interval $J_i$. If $B_{\mathbf{u}}$ were empty, then (35) would follow directly.

By Proposition 2, $P_f(G^c) < \varepsilon$ for all sufficiently large $n$. On the complementary event $G$, the Lebesgue measure of the set $B_{\mathbf{u}}$ of "bad" intervals $J_i$ is $< 2\delta$. Because the intervals $J_i$ *not* contained in $B_{\mathbf{u}}$ must have approximately the expected frequency $n|J_i|$ of data points, by (30), the number of data points in $B_{\mathbf{u}}$ cannot exceed $4\delta$, on the event $G$. Since $1 \geq B(k, l) \geq 4^{-k-l}$, it follows that the summands in (36) for which $J_i \subset B_{\mathbf{u}}$ cannot contribute more than $4\delta \log 4$ to the right-hand side. Assertion (35) now follows (with a larger value of $\delta$). □

COROLLARY 4. *For each $\delta > 0$, there exist $\varepsilon > 0$ and $K < \infty$ such that the following is true: If $K \leq m \leq \varepsilon n$ and $n$ is sufficiently large, then with $P_f$-probability at least $1 - \delta$,*

$$(37) \qquad Q_m(\{g : \|g - f\|_1 \geq \delta\} | (\mathbf{X}, \mathbf{Y})_n) < \delta$$

*and*

$$(38) \qquad |n^{-1} \log Z_m((\mathbf{X}, \mathbf{Y})_n) + H(f)| < \delta.$$

PROOF. For large $m$ (say, $m \geq K$), most split-point vectors $\mathbf{u}$ (as measured by the uniform distribution on $[0, 1]^m$) are such that all intervals $J_i(\mathbf{u})$ in the induced partition are short—this follows, for instance, from the Glivenko–Cantelli theorem—and so, by Lemma 1, $\|f - \bar{f}_{\mathbf{u}}\|_1$ is small. Thus, for any $\alpha, \varepsilon > 0$, there exists $K < \infty$ such that if $m \geq K$, then the set

$$B_m(\alpha) := \{\mathbf{u} \in [0, 1]^m : \|f - \bar{f}_{\mathbf{u}}\|_1 \geq \alpha\}$$

has Lebesgue measure $< \varepsilon$. Inequality (34) of Proposition 3 implies that, for each $\mathbf{u}$ in the complementary event $B_m^c(\alpha)$, the $Q_{\mathbf{u}}$-posterior distribution is concentrated on a small $L^1$-neighborhood of $f$, provided $\alpha$ is small. Thus, to prove (37), it must be shown that the contribution to the $Q_m$-posterior (14) from split-point vectors $\mathbf{u} \in B_m(\alpha)$ is negligible. For this, it suffices to show that the predictive probabilities $Z_{\mathbf{u}}((\mathbf{X}, \mathbf{Y})_n)$ are not larger for $\mathbf{u} \in B_m(\alpha)$ than for $\mathbf{u} \in B_m^c(\alpha)$.

By the entropy concavity property, $H(\bar{f}_{\mathbf{u}}) \geq H(f)$ for all split-point vectors $\mathbf{u}$, and for $\mathbf{u} \in B_m(\alpha)$,

$$H(\bar{f}) > H(f) + C\alpha,$$

where $2C = \max H''(p)$ [see (27)]. On the other hand, by the Entropy continuity property, if $\rho > 0$ is sufficiently small, then for all $\mathbf{u} \notin B_m(\rho)$,

$$|H(f) - H(\bar{f}_{\mathbf{u}})| < C\alpha/2.$$

By the second assertion (35) of Proposition 3, for all sufficiently large $n$ there is $P_f$-probability at least $1 - \delta$ that $n^{-1} \log Z_{\mathbf{u}}((\mathbf{X}, \mathbf{Y})_n)$ is within $\delta$ of $-H(\bar{f}_{\mathbf{u}})$ for all $\mathbf{u}$. Consequently, the primary contribution to the integral (14) must come from $\mathbf{u} \notin B_m(\alpha)$. This proves (37). Assertion (38) also follows, in view of the representation (15) for the predictive probability $Z_m((\mathbf{X}, \mathbf{Y})_n)$. $\square$

The exponential decay rate of the predictive probabilities in the beginning zone depends on whether or not the true regression function $f$ is a step function. If not, the decay is faster than in the middle zone; if so, the decay matches that in the middle zone, but the posterior concentrates in a neighborhood of $f$.

COROLLARY 5. *If the regression function $f$ is a step function with $k$ discontinuities in $(0, 1)$, then for each $m \geq k$ and all $\varepsilon > 0$, inequalities (37) and (38) hold with $P_f$-probability tending to 1 as the sample size $n \to \infty$. If $f$ is not a step function with fewer than $K + 1$ discontinuities, then there exists $\varepsilon > 0$ such that, with $P_f$-probability $\to 1$ as $n \to \infty$,*

$$(39) \qquad \max_{m \leq K} Z_m((\mathbf{X}, \mathbf{Y})_n) < \exp\{-nH(f) - n\varepsilon\}.$$

PROOF. If $f$ is not a step function with fewer than $K + 1$ discontinuities, then by the Entropy concavity property there exists $\varepsilon > 0$ so that $H(\bar{f}_{\mathbf{u}})$ is bounded above by $H(f) + \varepsilon$ for all split-point vectors $\mathbf{u}$ of length $m \leq K$. Hence, (39) follows from (35), by the same argument as in the proof of Corollary 4.

Suppose then that $f$ is a step function with $k$ discontinuities, that is, $f = \bar{f}_{\mathbf{u}^*}$ for some split-point vector $\mathbf{u}^*$ of length $k$. For any other split-point vector $\mathbf{u}$, the entropy $H(\bar{f}_{\mathbf{u}})$ cannot exceed $H(f)$ by the Entropy concavity property, and so (35) implies that, for any $m$, the exponential decay rate of the predictive probability $Z_m((\mathbf{X}, \mathbf{Y})_n)$ as $n \to \infty$ cannot exceed $-H(f)$. But since $f$ is a step function with $k$ discontinuities, any open $L^1$ neighborhood of $f$ has positive $\pi_m$-probability; consequently, by entropy continuity and (35), the exponential decay rate of $Z_m((\mathbf{X}, \mathbf{Y}_n))$ in $n$ must be at least $-H(f)$. Thus, (38) holds with $P_f$-probability $\to 1$ as $n \to \infty$. Finally, (39) follows by the same argument as in the proof of Corollary 4. □

Corollaries 5 and 4 imply that, with $P_f$-probability $\to 1$ as $n \to \infty$, the $Q^\nu$-posterior in the beginning and middle zones concentrates near $f$, and that the total posterior mass in the beginning and middle zones decays at the exponential rate $H(f)$ as $n \to \infty$. Thus, to complete the proof of Theorem 1, it suffices to show that the posterior mass in the end zone $m \geq \delta n$ decays at an exponential rate $> H(f)$. This will be the agenda for the remainder of the article: see Proposition 6 below.

**4. The end zone.** For the Diaconis–Freedman priors, the log-predictive probabilities simplify neatly as sums of independent random variables, and so their asymptotic behavior drops out easily from the usual WLLN. No such simplification is possible in our case: the integral in (15) does not admit further reduction. Thus, the analysis of the posterior in the end zone will necessarily be somewhat more roundabout than in the Diaconis–Freedman case. The main objective is the following.

PROPOSITION 6. *For any Borel measurable regression function $f \not\equiv 1/2$ and all $\varepsilon > 0$, there exist constants $\delta = \delta(\varepsilon, f) > 0$ such that*

$$(40) \qquad \lim_{n \to \infty} P_f \left\{ \sup_{m \geq \varepsilon n} \log Z_m((\mathbf{X}, \mathbf{Y})_n) \geq n(-H(f) - \delta) \right\} = 0.$$

Given this result, the consistency theorem follows.

PROOF OF THEOREM 1. Proposition 6 implies that, for any $\varepsilon > 0$, with high $P_f$-probability the posterior mass (un-normalized) in the region $m \geq \varepsilon n$ is less than $\exp\{-nH(f) - n\delta\}$. Corollaries 4 and 5 imply that there exists $\varepsilon' > 0$ so that, with $P_f$-probability tending to one as $n \to \infty$, the posterior mass in the region $m \leq \varepsilon' n$ is at least $\exp\{-nH(f) - n\delta/2\}$. Consequently, for any $\varepsilon > 0$, the posterior mass will, for large $n$, with high $P_f$-probability be almost entirely concentrated on step functions with $m \leq \varepsilon n$ discontinuities.

Assertion (37) of Corollary 4, together with Corollary 5, implies that for some $\varepsilon'' > 0$, most of the posterior mass in the region $m \leq \varepsilon'' n$ will, with high $P_f$-probability be concentrated on step functions near $f$. Since by the preceding paragraph nearly all of the posterior mass will eventually be concentrated in the region $m \leq \varepsilon'' n$, the result (1) follows. □

To prove Proposition 6, we will show in Proposition 11 below that the predictive probabilities (after suitable "Poissonization") decay exponentially in $n$ at a precise rate, depending on $\alpha > 0$, for $m/n \to \alpha > 0$.

4.1. *Preliminaries*: *Comparison and Poissonization*.  Comparison arguments will be based on the following simple observation.

LEMMA 7.  *Adding more data points* $(x_i, y_i)$ *to the sample* $(\mathbf{x}, \mathbf{y})$ *decreases the value of* $Z_m(\mathbf{x}, \mathbf{y})$.

PROOF.  For each fixed pair $\mathbf{u}, \mathbf{w} \in [0, 1]^m$, adding data points to the sample increases at least some of the cell counts $N_i^S, N_i^F$, and therefore, decreases the integrand in (12).  □

Two Poissonizations will be used, one for the data sample, the other for the sample of split points. Let $\Lambda(t)$ be a standard Poisson counting process of intensity 1, independent of the data stream $(\mathbf{X}, \mathbf{Y})$. Replacing the sample $(\mathbf{X}, \mathbf{Y})_n$ of fixed size $n$ by a sample $(\mathbf{X}, \mathbf{Y})_{\Lambda(n)}$ of size $\Lambda(n)$ has the effect of making the success/failure counts in disjoint intervals independent random variables with Poisson distributions.

LEMMA 8.  *For each* $\varepsilon > 0$, *the probability that*

$$(41) \qquad Z_m((\mathbf{X}, \mathbf{Y})_{\Lambda(n-\varepsilon n)}) \le Z_m((\mathbf{X}, \mathbf{X})_n) \le Z_m((\mathbf{X}, \mathbf{Y})_{\Lambda(n+\varepsilon n)})$$

*for all* $m$ *approaches* 1 *as* $n \to \infty$.

PROOF.  For any $\varepsilon > 0$, $P\{\Lambda(n - \varepsilon n) \le n \le \Lambda((1 + \varepsilon)n)\} \to 1$ as $n \to \infty$, by the weak law of large numbers. On this event, inequality (41) must hold by Lemma 7.  □

The second Poissonization, for the split point vector, involves mixing the priors $\pi_m$ according to a Poisson hyperprior. For any $\lambda > 0$, let $\pi_\lambda^*$ be the Poisson-$\lambda$ mixture of the priors $\pi_m$, and let $Q_\lambda^*$ be the corresponding induced measure on data sequences (equivalently, $Q_\lambda^*$ is the Poisson-$\lambda$ mixture of the measures $Q_m$). Then the $Q_\lambda^*$-predictive probability for a data set $(\mathbf{x}, \mathbf{y})$ is given by

$$(42) \qquad\qquad Z_\lambda^*(\mathbf{x}, \mathbf{y}) := \sum_{k=0}^{\infty} \frac{\lambda^k e^{-\lambda}}{k!} Z_k(\mathbf{x}, \mathbf{y}).$$

The effect of Poissonization on the number of split points is a bit more subtle than the effect on data, because there is no simple a priori relation between neighboring predictive probabilities $Z_m(\mathbf{x}, \mathbf{y})$ and $Z_{m+1}(\mathbf{x}, \mathbf{y})$. However, because the Poisson distribution with mean $\alpha n$ assigns mass at least $C/\sqrt{n}$ to the value $[\alpha n]$ (by the Local CLT), where $C = C(\alpha) > 0$ is continuous in $\alpha$, the following is obviously true.

LEMMA 9. *For each $\varepsilon > 0$, there exists $C < \infty$ such that, for all $\alpha \in (\varepsilon, \varepsilon^{-1})$,*

$$(43) \qquad Z_{[\alpha n]}((\mathbf{X}, \mathbf{Y})_{\Lambda(n)}) \leq C\sqrt{n} Z^*_{\alpha n}((\mathbf{X}, \mathbf{Y})_{\Lambda(n)}).$$

*Consequently, to prove Proposition 6 it suffices to prove that (40) holds when $Z_m((\mathbf{X}, \mathbf{Y})_n)$ is replaced by $Z^*_m((\mathbf{X}, \mathbf{Y})_n)$.*

Whereas it is difficult to compare neighboring predictive probabilities $Z_m(\mathbf{x}, \mathbf{y})$ and $Z_{m+1}(\mathbf{x}, \mathbf{y})$, it is quite easy to compare *Poissonized* predictive probabilities $Z^*_\lambda(\mathbf{x}, \mathbf{y})$ and $Z^*_\mu(\mathbf{x}, \mathbf{y})$ for neighboring intensities $\mu, \lambda$.

LEMMA 10. *For each $\varepsilon > 0$ and $A < \infty$, there exists $\delta > 0$ such that if $\mu, \lambda \leq A$ and $|\mu - \lambda| \leq \delta$, then for all $n \geq 1$ and all data sets $(\mathbf{x}, \mathbf{y})$ of size $n$,*

$$(44) \qquad \frac{Z^*_{\mu n}(\mathbf{x}, \mathbf{y})}{Z^*_{\lambda n}(\mathbf{x}, \mathbf{y})} \leq e^{n\varepsilon}.$$

PROOF. Inequality (22) implies that $Z^*_m(\mathbf{x}, \mathbf{y}) \geq 4^{-n}$. Chernoff's large deviation inequality implies that if $M$ has the Poisson distribution with mean $\lambda \leq An$, then

$$P\{M \geq \kappa n\} \leq e^{-\gamma n},$$

where $\gamma \to \infty$ as $\kappa \to \infty$. Since $Z_k(\mathbf{x}, \mathbf{y}) \leq 1$, it follows that the contribution to the sum (42) from terms indexed by $k \geq \kappa n$ is of smaller exponential order of magnitude than that from terms indexed by $k < \kappa n$, provided $\gamma > \log 4$.

Consider the Poisson distributions with means $\mu n, \lambda n \leq \kappa n$: these are mutually absolutely continuous, and the likelihood ratio at the integer value $k$ is

$$(\mu/\lambda)^k e^{n\lambda - n\mu}.$$

If $k \leq \kappa n$ and $|\mu - \lambda|$ is sufficiently small, then this likelihood ratio is less than $e^{n\varepsilon}$. By the result of the preceding paragraph, only values of $k \leq \kappa n$ contribute substantially to the expectations; thus, the assertion follows. □

In some of the arguments to follow, an alternative representation of these Poissonized predictive probabilities as a conditional expectation will be useful. Assume that on the underlying probability space $(\Omega, \mathcal{F}, P_f)$ are defined i.i.d. uniform-$[0,1]$ r.v.s $U_n$ and independent Poisson processes $\Lambda, M$, all jointly independent of the data stream. Then

$$(45) \qquad Z^*_\lambda((\mathbf{X}, \mathbf{Y^f})_{\Lambda(n)}) = E_f(\beta | (\mathbf{X}, \mathbf{Y})_{\Lambda(n)}),$$

where

$$\beta = \beta(\mathbf{U}_{M(\lambda)}; (\mathbf{X}, \mathbf{Y^f})_{\Lambda(n)}) := \prod_{i=0}^{M(\lambda)} B(N_i^S, N_i^F) \tag{46}$$

and $N_i^S, N_i^F$ are the success/failure cell counts for the data $(\mathbf{X}, \mathbf{Y})_{\Lambda(n)}$ relative to the partition induced by the split point sample $\mathbf{U}_{M(\lambda)}$.

4.2. *Exponential decay.* The asymptotic behavior of the doubly Poissonized predictive probabilities is spelled out in the following proposition, whose proof will be the goal of Sections 4.4–4.6 and Section 5 below.

PROPOSITION 11. *For each Borel measurable regression function $f$ and each $\alpha > 0$, there exists a constant $\psi_f(\alpha)$ such that, as $n \to \infty$,*

$$n^{-1} \log Z^*_{\alpha n}(((\mathbf{X}, \mathbf{Y})_{\Lambda(n)})) \xrightarrow{P_f} \psi_f(\alpha). \tag{47}$$

*The function $\psi_f(\alpha)$ satisfies*

$$\psi_f(\alpha) = \int_0^1 \psi(f(x), \alpha) \, dx, \tag{48}$$

*where $\psi(p, \alpha) = \psi_p(\alpha)$ is the corresponding limit for the constant regression function $f \equiv p$. The function $\psi(p, \alpha)$ is jointly continuous in $p, \alpha$ and satisfies*

$$\lim_{\alpha \to \infty} \max_{p \in [0,1]} |\psi_p(\alpha) + \log 2| = 0 \tag{49}$$

*and*

$$\psi_p(\alpha) < -H(p). \tag{50}$$

Note that the entropy inequality (50) extends to all regression functions $f$: that is, $p$ may be replaced by $f$ on both sides of (50). This follows from the integral formulas that define $\psi_f(\alpha)$ and $H(f)$. The fact that this inequality is strict is crucially important to the consistency theorem. It will also require a rather elaborate argument: see Section 5 below.

The case $f \equiv p$, where the regression function is constant, will prove to be the crucial one. In this case, the existence of the limit (47) is somewhat reminiscent of the existence of "thermodynamic limits" in formal statistical mechanics (see [21], Chapter 3). Unfortunately, Proposition 11 cannot be reduced to the results of [21], as follows: (i) the data sequence enters conditionally (thus functioning as a "random environment"); and, more importantly, (ii) the hypothesis of "tempered interaction" needed in [21] cannot be verified here. The limit (47) is also related to the "conditional LDP" of Chi [2], but again cannot be deduced from the results of that paper, because the log-predictive probability cannot be expressed as a continuous functional of the empirical distribution of split point/data point pairs.

4.3. *Proof of Proposition* 6. Before proceeding with the somewhat arduous proof of Proposition 11, we show how to complete the proof of Proposition 6. In the process, we shall establish the asymptotic behavior (49) of the rate function.

LEMMA 12. *For every $\delta > 0$, there exists $\alpha_\delta < \infty$ such that*

$$\lim_{n\to\infty} P\left\{\sup_{\alpha \geq \alpha_\delta} \sup_{\mathbf{y}} Z^*_{\alpha n}((\mathbf{X}, \mathbf{y})_{\Lambda(n)}) \geq 2^{-n+n\delta}\right\} = 0 \tag{51}$$

*and*

$$\lim_{n\to\infty} P\left\{\inf_{\alpha \geq \alpha_\delta} \inf_{\mathbf{y}} Z^*_{\alpha n}((\mathbf{X}, \mathbf{y})_{\Lambda(n)}) \leq 2^{-n-n\delta}\right\} = 0. \tag{52}$$

*Here $\sup_{\mathbf{y}}$ and $\inf_{\mathbf{y}}$ are taken over all assignments of 0s and 1s to the response variables $y_1, y_2, \ldots, y_{\Lambda(n)}$. Similarly,*

$$\lim_{n\to\infty} P\left\{\sup_{m \geq \alpha_\delta n} \sup_{\mathbf{y}} Z_m((\mathbf{X}, \mathbf{y})_n) \geq 2^{-n+n\delta}\right\} = 0 \tag{53}$$

*and*

$$\lim_{n\to\infty} P\left\{\inf_{m \geq \alpha_\delta n} \inf_{\mathbf{y}} Z_m((\mathbf{X}, \mathbf{y})_n) \leq 2^{-n-n\delta}\right\} = 0. \tag{54}$$

Given the convergence (47), the following is now immediate from (51)–(52).

COROLLARY 13. *For every regression function $f$,*

$$\lim_{\alpha\to\infty} \psi_f(\alpha) = -\log 2. \tag{55}$$

PROOF OF LEMMA 12. We shall prove only (51)–(52); the other two assertions may be proved by similar arguments. Let $\xi_1, \xi_2, \ldots, \xi_{\Lambda(n)+1}$ be the spacings between successive order statistics of the covariates $X_1, X_2, \ldots, X_{\Lambda(n)}$. For each pair of positive reals $\varepsilon, \delta > 0$, let $G = G_{\delta,\varepsilon}$ be the event that at least $(1 - \delta)n$ of the spacings $\xi_j$ are larger than $\varepsilon/n$. Call these spacings "fat." Since the spacings are independent exponentials with mean $1/n$, the Glivenko–Cantelli theorem implies that there exist $\delta = \delta(\varepsilon) \to 0$ as $\varepsilon \to 0$ such that

$$\lim_{n\to\infty} P(G_{\delta,\varepsilon} \cap \{|\Lambda(n) - n| < \varepsilon n\}) = 1.$$

By elementary large deviations estimates for the Poisson process, given $G$, the probability that a random sample of $M(\alpha n)$ split points is such that more than $(1 - 2\delta)n$ of the fat spacings contain no split points is less than $\exp\{-n\gamma\}$, where $\gamma = \gamma(\alpha, \varepsilon, \delta) \to \infty$ as $\alpha \to \infty$. But on the complement of this event, at least $(1 - 4\delta)n$ of the intervals induced by the split points have *exactly* one data point. Thus, on the event $G \cap \{|\Lambda(n) - n| < \varepsilon n\}$,

$$2^{-n+4n\delta}4^{-4\delta n - \varepsilon n} \leq \prod_i B(N_i^S, N_i^F) \leq 2^{-n+4\delta n}.$$

Observe that these inequalities hold regardless of the assignment $\mathbf{y}$ of values to the response variables. Thus, taking conditional expectations (45) given the data $(\mathbf{X}, \mathbf{Y})_{\Lambda(n)}$, we obtain

$$(56) \quad (1 - e^{-n\gamma})2^{-n-4n\delta-2\delta n} \leq Z_{\alpha n}^*(\mathbf{X}_{\Lambda(n)}, \mathbf{Y}_{\Lambda(n)}) \leq 2^{-n+4\delta n} + e^{-n\gamma}.$$

Since $\gamma$ can be made arbitrarily large by making $\alpha$ large, assertions (51) and (52) follow. $\square$

PROOF OF PROPOSITION 6. Since $H(f) < \log 2$ for every regression function $f \not\equiv 1/2$, Lemma 12 implies that, to prove (40), it suffices to replace the supremum over $m \geq \varepsilon n$ by the supremum over $m \in [\varepsilon n, \varepsilon^{-1} n]$. Now for $m$ in this range, the bound (43) is available; since $\log n$ is negligible compared to $n$, (43) implies that

$$\sup_{\varepsilon n \leq m \leq \varepsilon^{-1} n} \log Z_m((\mathbf{X}, \mathbf{Y})_n)$$

may be replaced by

$$\sup_{\varepsilon \leq \alpha \leq \varepsilon^{-1}} \log Z_{\alpha n}^*((\mathbf{X}, \mathbf{Y})_{\Lambda(n)})$$

in (40). Lemma 10 implies that this last supremum may be replaced by a maximum over a finite set of values $\alpha$, and now (40) follows from assertions (47), (48) and (50) of Proposition 11. $\square$

4.4. *Constant regression functions.* The simplest route to the convergence (47) is via subadditivity (more precisely, approximate subadditivity) arguments. Assume that $f \equiv p$ is constant, and that the constant $p \neq 0, 1$. Recall (Section 2.4) that, in this case, the integral (15) defining the predictive probability almost factors perfectly into the product of two integrals, one over the data and split points in $[0, 1/2]$, the other over $(1/2, 1]$, of the same form [but on a different scale—see (18)]. Unfortunately, this factorization is not exact, as the partition of the unit interval induced by the split points $u_i$ includes an interval that straddles the demarcation point $1/2$. However, the error can be controlled, and so the convergence (47) can be deduced from a subadditive WLLN (Proposition A.1 of the Appendix A). The next lemma shows that the hypotheses of Proposition A.1 are met.

LEMMA 14. *Fix $\alpha > 0$, and write $\zeta_n = \log Z_{\alpha n}^*(((\mathbf{X}, \mathbf{Y})_{\Lambda(n)}))$. For each pair $m, n \in \mathbb{N}$ of positive integers, there exist random variables $\zeta'_{m,m+n}$, $\zeta''_{n,m+n}$ and $R_{m,n}$ such that:*

(a) $\zeta'_{m,m+n}, \zeta''_{n,m+n}$ *are independent;*
(b) $\zeta_m$ *and* $\zeta'_{m,m+n}$ *have the same law;*
(c) $\zeta_n$ *and* $\zeta''_{n,m+n}$ *have the same law;*
(d) *the random variables* $\{R_{m,n}\}_{m,n \geq 1}$ *are identically distributed;*
(e) $E|R_{1,1}| < \infty$; *and*
(f) *for all integers* $m, n \geq 1$,

$$(57) \qquad \zeta_{m+n} \geq \zeta'_{m,m+n} + \zeta''_{n,m+n} + R_{m,n}.$$

Proposition A.1 of the Appendix A now implies the following.

COROLLARY 15. *For some constant* $\psi_p(\alpha)$,

$$(58) \qquad n^{-1}\zeta_n \xrightarrow{L^1} \psi_p(\alpha).$$

PROOF OF LEMMA 14. The construction requires auxiliary randomization, and so we assume that independent copies of the data sequence and the various Poisson processes are available. Consider the expectation (45) that defines the Poissonized predictive probability $\exp\{\zeta_{m+n}\}$. This expectation extends over all samples of split points of size $M(\alpha m + \alpha n)$. Since the integrand is positive, the expectation exceeds its restriction to the event $G$ that there are split points in both of the intervals $[b - (m+n)^{-1}, b]$ and $(b, b + (m+n)^{-1}]$, where $b := m/(m+n)$. Note that this event has probability

$$\delta = \delta(\alpha) := (1 - e^{-\alpha})^2.$$

Denote by $U', U''$ the split points nearest the demarcation point $b$ to its left and right, respectively. The product $\beta = \prod B(N_i^S, N_i^F)$ may be factored into three parts, consisting of terms indexed by intervals $J_i$ contained in $[0, U']$, intervals contained in $(U'', 1]$ and the single interval $J_* = (U', U'')$ that straddles the point $b$. Conditional on the values of $U', U''$, the three products are independent: by the scaling relation (18), the first two have the same distributions as the products $\beta$ occurring as integrands in the expectations defining

$$\exp\{\zeta_{U'(m+n)}\} \quad \text{and} \quad \exp\{\zeta_{(1-U'')(m+n)}\},$$

respectively; and the third is just

$$B(N_*^S, N_*^F) \geq 2^{-N_*}/(N_* + 1) \geq 4^{-N_*},$$

where $N_*^S$ and $N_*^F$ are the numbers of successes and failures in the interval $J_*$, and $N_* = N_*^S + N_*^F$. Note that on $G$ the random variable $N_*$ is dominated by a Poisson random variable $N_{**}$ with mean 2, since the length of $J_*$ is less than $2/(m+n)$ (this requires auxiliary randomization). Now extend each of the first two products in the following manner: throw new, independent data samples and split points into the intervals $(U', b]$ and $(b, U'')$; remove the split points $U', U''$ and place a split at $b$; then recompute the partitions and replace the affected terms $B(N_i^S, N_i^F)$ by the new values. Note that this cannot increase the value of either product; moreover, the products remain conditionally independent given $U', U''$. Most importantly, the conditional expectations of these products (given the data and the values of $U', U''$) have the same distributions as $\exp\{\zeta_m\}$ and $\exp\{\zeta_n\}$, respectively. Thus,

$$\exp\{\zeta_{m+n}\} \geq \delta \exp\{\zeta_m'\} \exp\{\zeta_n''\} 4^{-N_{**}},$$

where $\zeta_m'$ and $\zeta_n''$ are independent, with the same distributions as $\zeta_m$ and $\zeta_n$, respectively, and $N_{**}$ is Poisson with mean 2.   □

REMARK.   There is a similar (and in some respects simpler) approximate *subaddivitivity* relation among the distributions of the random variables $\zeta_n$: For each pair $m, n \geq 1$ of positive integers, there exist independent random variables $\xi_{m,m+n}', \xi_{n,m+n}''$ whose distributions are the same as those of $\zeta_m, \zeta_n$, respectively, such that

$$(59) \qquad \zeta_{m+n} \leq \xi_{m,m+n}' + \xi_{n,m+n}'' + \log \Lambda(m+n).$$

Corollary 15 can also be deduced from (59), but this requires a more sophisticated subadditive LLN than is proved in Appendix A, because the remainders $\log \Lambda(m+n)$ are not uniformly $L^1$ bounded, as they are in (57). This approach has the advantage that it leads to a proof that the convergence (58) holds almost surely.

PROOF OF (59).   Consider the effect on the integral (15) of adding a split point at $b = m/(m+n)$: This breaks one of intervals $J_i$ into two, leaving all of the others unchanged, and so the effect on the integrand in (15) is that one of the factors $B(N_i^S, N_i^F)$ is replaced by a product of two factors $B(N_L^S, N_L^F)B(N_R^S, N_R^F)$. By (23), the multiplicative error in this replacement is bounded above by $\Lambda(m+n)$. After the replacement, the factors in the integrand $\beta = \prod B(N_i^S, B_i^F)$ may be partitioned neatly into those indexed by intervals left of $b$ and those indexed by intervals right of $b$: thus,

$$\beta = \beta'\beta'',$$

where $\beta, \beta''$ are independent and have the same distributions as the products $\beta$ occurring as integrands in the expectations defining $\exp\{\zeta_m\}$ and $\exp\{\zeta_n\}$, respectively. Thus,

$$\exp\{\zeta_{m+n}\} \leq \exp\{\xi_m\} \exp\{\xi_n''\} \Lambda(m+n),$$

where $\xi_m' = \xi_{m,m+n}'$ and $\xi_n'' = \xi_{n,m+n}''$ are independent and distributed as $\zeta_m$ and $\zeta_n$, respectively.   □

4.5. *Piecewise constant regression functions.* The next step is to extend the convergence (47) to piecewise constant regression functions $f$. For ease of exposition, we shall restrict attention to step functions with a single discontinuity in $(0, 1)$; the general case involves no new ideas. Thus, assume that

$$f(x) = p_L \qquad \text{for } x \leq b,$$
$$f(x) = p_R \qquad \text{for } x > b,$$

and

$$p_L \neq p_R.$$

Fix $\alpha > 0$, and set

$$Z_n^* := Z_{\alpha n}^*((\mathbf{X}, \mathbf{Y})_{\Lambda(n)}). \tag{60}$$

LEMMA 16. *With $P_f$-probability approaching one as $n \to \infty$,*

$$Z_n^* \geq Z_n' Z_n'' / n^2 \tag{61}$$

*and*

$$Z_n^* \leq 2n Z_n' Z_n'', \tag{62}$$

*where, for each $n$, the random variables $Z_n', Z_n''$ are independent, with the same distributions as*

$$Z_n' \overset{\mathcal{L}}{=} Z_{\alpha n b}^*((\mathbf{X}, \mathbf{Y})_{\Lambda(bn)}) \qquad \text{under } P_{p_L};$$
$$Z_n'' \overset{\mathcal{L}}{=} Z_{\alpha n - \alpha n b}^*((\mathbf{X}, \mathbf{Y})_{\Lambda(n-nb)}) \qquad \text{under } P_{p_R}. \tag{63}$$

PROOF. Consider the effect on $Z_n^*$ of placing an additional split point at $b$: this would divide the interval straddling $b$ into two nonoverlapping intervals $L, R$ (for "left" and "right"), and so in the integrand $\beta := \prod B(N_i^S, N_i^F)$ the single factor $B(N_*^S, N_*^F)$ representing the interval straddling $b$ would be replaced by a product of two factors $B(N_L^S, N_L^F)$ and $B(N_R^S, N_R^F)$. As in the proof of the subadditivity inequality (59) in Section 4.4, the factors of this modified product separate into those indexed by subintervals of $[0, b]$ and those indexed by subintervals of $[b, 1]$; thus, the modified product has the form $\beta' \beta''$, where $\beta'$ and $\beta''$ are the products of the factors indexed by intervals to the left and right, respectively, of $b$. Denote by $Z_n'$ and $Z_n''$ the conditional expectations of $\beta'$ and $\beta''$ (given the data). These are independent random variables, and by the scaling relation (18), their distributions satisfy (63). By inequality (23), the multiplicative error in making the replacement is at most $\Lambda(n)$; since the event $\Lambda(n) \geq 2n$ has probability tending to 0 as $n \to \infty$, inequality (62) follows.

The reverse inequality (61) follows by a related argument. Let $G$ be the event that the data sample $(\mathbf{X}, \mathbf{Y})_{\Lambda(n)}$ contains no points with covariate $X_i \in [b, b + n^{-2}]$. Since the covariates are generated by a Poisson point process with intensity $n$, the probability of $G^c$ is approximately $n^{-1}$. Consider the integral (over all samples of split points) that defines $Z_n^*$: this integral exceeds its restriction to the event $A$ that there is a split point in $[b, b + n^{-2}]$. The conditional probability of $A$ (given the data) is approximately $\alpha n^{-1}$, and thus larger than $n^{-2}$ for large $n$. On the event $G \cap A$,

$$\beta = \beta' \beta''$$

holds exactly, as the split point in $[b, b + n^{-2}]$ produces exactly the same bins as if the split point were placed at $b$. Moreover, conditioning on the event $A$ does not affect the joint distribution (conditional on the data) of $\beta', \beta''$ when $G$ holds. Thus, the conditional expectation of the product, given $A$ and the data, equals $Z_n' Z_n''$ on the event $G$.    $\square$

Taking $n$th roots on each side of (61) and appealing to Corollary 15 now yields the following.

COROLLARY 17. *If the regression function is piecewise constant, with only finitely many discontinuities, then the convergence (47) holds.*

4.6. *Thinning.* Extension of the preceding corollary to arbitrary Borel measurable regression functions will be based on *thinning* arguments. Recall that if points of a Poisson point process of intensity $\lambda(x)$ are randomly removed with location-dependent probability $\varrho(x)$, then the resulting "thinned" point process is again Poisson, with intensity $\lambda(x) - \varrho(x)\lambda(x)$. This principle may be applied to both the success $(y = 1)$ and failure $(y = 0)$ point processes in a Poissonized data sample. Because thinning at location-dependent rates may change the distribution of the covariates, it will be necessary to deal with data sequences with nonuniform covariate distribution. Thus, let $(\mathbf{X}, \mathbf{Y})$ be a data sample of random size with Poisson-$\lambda$ distribution under the measure $P_{f,F}$ (here $f$ is the regression function, $F$ is the distribution of the covariate sequence $X_j$). If successes $(x, y = 1)$ are removed from the sample with probability $\varrho_1(x)$ and failures $(x, y = 0)$ are removed with probability $\varrho_0(x)$, then the resulting sample will be a data sample of random size with the Poisson-$\mu$ distribution from $P_{g,G}$, where the mean $\mu$, the regression function $g$ and the covariate distribution $G$ satisfy

$$(64) \qquad \begin{aligned} \mu g(x) G(dx) &= (1 - \varrho_1(x))\lambda f(x) F(dx) \quad \text{and} \\ \mu(1 - g(x)) G(dx) &= (1 - \varrho_0(x))\lambda(1 - f(x)) F(dx). \end{aligned}$$

By the monotonicity principle (Lemma 7), the predictive probability of the thinned sample will be no smaller than that of the original sample. Thus, thinning allows comparison of predictive probabilities for data generated by two different measures $P_{f,F}$ and $P_{g,G}$. The first and easiest consequence is the continuity of the rate function.

LEMMA 18.  *The rate function* $\psi_p(\alpha)$ *is jointly continuous in* $p, \alpha$.

PROOF.  Corollary 15 and Lemma 10 imply that the functions $\alpha \mapsto \psi_p(\alpha)$ are uniformly continuous in $\alpha$. Continuity in $p$ and joint continuity in $p, \alpha$ are now obtained by thinning. Let $(\mathbf{X}, \mathbf{Y})$ be a random sample of size $\Lambda(n) \sim$ Poisson-$n$ from a data stream distributed according to $P_p$ (i.e., $f \equiv p$ and $F$ is the uniform-$[0, 1]$ distribution). Let $(\mathbf{X}, \mathbf{Y})'$ be the sample obtained by randomly removing failures from the sample $(\mathbf{X}, \mathbf{Y})$, with probability $\varepsilon$. Then $(\mathbf{X}, \mathbf{Y})'$ has the same distribution as a random sample of size $\Lambda(n - \varepsilon q n)$ (here $q = 1 - p$) from a data stream distributed according to $P_{p'}$, where $p' = p/(1 - \varepsilon q)$. By the monotonicity principle (Lemma 7),

$$Z^*_{\alpha n}((\mathbf{X}, \mathbf{Y})) \leq Z^*_{\alpha n}((\mathbf{X}, \mathbf{Y})').$$

Taking $n$th roots and appealing to Corollary 15 shows that

$$\psi(p, \alpha) \leq (1 - \varepsilon q)^{-1} \psi(p/(1 - \varepsilon q), \alpha/(1 - \varepsilon q)).$$

A similar inequality in the opposite direction can be obtained by reversing the roles of $p$ and $p/(1 - \varepsilon q)$. The continuity in $p$ of $\psi(p, \alpha)$ now follows from the continuity in $\alpha$, and the joint continuity follows from the uniform continuity in $\alpha$.  □

PROPOSITION 19.  *The convergence* (47) *holds for every Borel measurable regression function* $f$.

PROOF.  By Corollary 17 above, the convergence holds for all piecewise constant regression functions with only finitely many discontinuities. The general case will be deduced from this by another thinning argument.

If $f : 0, 1 \to [0, 1]$ is measurable, then for each $\varepsilon > 0$, there exists a piecewise constant $g : [0, 1] \to [0, 1]$ (with only finitely many discontinuities) such that $\|f - g\|_1 < \varepsilon$. If $\varepsilon$ is small, then $|f - g|$ must be small except on a set $B$ of small Lebesgue measure; moreover, $g$ may be chosen so that $g = 1$ wherever $f$ is near 1, and $g = 0$ wherever $f$ is near 0 (except on $B$). For such choices of $\varepsilon$ and $g$, there will exist removal rate functions $\varrho_0(x)$ and $\varrho_1(x)$ so that equation (64) holds with $F =$ the uniform distribution on $[0, 1]$, $G =$ the uniform distribution on $[0, 1] - B$, and

$$\left| \frac{\lambda}{\mu} - 1 \right| < \delta(\varepsilon)$$

for some constants $\delta(\varepsilon) \to 0$ as $\varepsilon \to 0$. (Note: Requiring $G$ to be the uniform distribution on $[0, 1] - B$ forces complete thinning in $B$, i.e., $\varrho_0 = \varrho_1 = 1$ in $B$.) Thus, a Poissonized data sample distributed according to $P_f$ may be thinned so as to yield a Poissonized data sample distributed according to $P_{g,G}$ in such a way that the overall thinning rate is arbitrarily small. It follows, by the monotonicity principle, that the Poissonized predictive probabilities for data distributed according to $P_f$ are majorized by those for data distributed according to $P_{g,G}$, with a slightly smaller rate.

Now consider data $(\mathbf{X}, \mathbf{Y})$ distributed according to $P_{g,G}$: Since $g$ is piecewise constant and $G$ is a uniform distribution, the transformed data $(G\mathbf{X}, \mathbf{Y})$ will be distributed as $P_h$, where $h$ is again piecewise constant. Moreover, since the removed set $B$ has small Lebesgue measure, the function $h$ is close to the function $g$ in the Skorohod topology, and so by Lemma 18, $\psi_h \approx \psi_g \approx \psi_f$. Because the convergence (47) has been established for piecewise constant regression functions $h$, it now follows from the monotonicity principle that

$$P_f\{n^{-1} \log Z_{\alpha n}^*(((\mathbf{X}, \mathbf{Y})_{\Lambda(n)})) > \psi_f(\alpha) + \delta\} \longrightarrow 0$$

for every $\delta > 0$. This proves the upper (and for us, the more important) half of (47). The lower half may be proved by a similar thinning argument in the reverse direction. $\square$

**5. Proof of the entropy inequality (50).** This requires a change of perspective. Up to now, we have taken the point of view that the covariates $X_j$ and the split points $U_i$ are generated by Poisson point processes in the unit interval of intensities $n$ and $\alpha n$, respectively. However, the transformation formula (18) implies that the predictive probabilities and hence also their Poissonized versions, are unchanged if the covariates and the split points are rescaled by a common factor $n$. The rescaled covariates $\hat{X}_j := X_j/n$ and split points $\hat{U}_i := U_i/n$ are then generated by Poisson point processes of intensities 1 and $\alpha$ on the interval $[0, n]$. Consequently, versions of all the random variables $Z_{[\alpha n]}^*((\mathbf{X}, \mathbf{Y})_{\Lambda(n)})$ may be constructed from two independent Poisson processes of intensities 1 and $\alpha$ on the whole real line. The advantage of this new point of view is the possibility of deducing the large-$n$ asymptotics from the Ergodic theorem.

5.1. *Reformulation of the inequality.* To avoid cluttered notation, we shall henceforth drop the hats from the rescaled covariates and split points. Thus, assume that under both $P = P_p$ and $Q$,

$$\cdots < X_{-1} < X_0 < 0 < X_1 < \cdots$$

and

$$\cdots < U_{-1} < U_0 < 0 < U_1 < \cdots$$

are the points of independent Poisson point processes $\mathbf{X}$ and $\mathbf{U}$ of intensities 1 and $\alpha$, respectively, and let $\{W_i\}_{i \in \mathbb{Z}}$ be a stream of uniform-$[0,1]$ random variables independent of the point processes $\mathbf{X}, \mathbf{U}$. Denote by $N(t)$ the number of occurrences in the Poisson point process $\mathbf{X}$ during the interval $[0,t]$, and set $J_i = (U_i, U_{i+1})$. Let $\{Y_i\}_{i \in \mathbb{Z}}$ be Bernoulli r.v.s distributed according to the following laws:

(A) Under $P$, the random variables $Y_j$ are i.i.d. Bernoulli-$p$, jointly independent of the Poisson point processes $\mathbf{U}, \mathbf{X}$.

(B) Under $Q$, the random variables $Y_j$ are conditionally independent, given $\mathbf{X}, \mathbf{U}, \mathbf{W}$, with conditional distributions $Y_j \sim$ Bernoulli-$W_i$, where $i$ is the index of the interval $J_i$ containing $X_j$.

Under $Q$, the sequence $\{Y_n\}_{n \in \mathbb{Z}}$ is an ergodic, stationary sequence; for reasons that we shall explain below, we shall refer to this process as the *rechargeable Pólya urn*. The distribution of $(\mathbf{X}, \mathbf{Y}) \cap [0,t]$ under $Q$ is, after rescaling of the covariates by the factor $t$, the same as that of a data sample of random size $\Lambda(t)$ under the Poisson mixture $Q^*_{\alpha t}$ defined in Section 4.1 above.

For (extended) integers $-\infty \leq m \leq n$, define $\sigma$-algebras

$$\mathcal{F}^{X,Y}_{m,n} = \sigma(\{X_j, Y_j\}_{m \leq j \leq n}) \quad \text{and} \quad \mathcal{F}^Y_{m,n} = \sigma(\{Y_j\}_{m \leq j \leq n}).$$

If $m, n$ are both finite, then the restrictions of the measures $P, Q$ to $\mathcal{F}^{X,Y}_{m,n}$ (and therefore also to $\mathcal{F}^Y_{m,n}$) are mutually absolutely continuous. The Radon–Nikodym derivative on the smaller $\sigma$-algebra $\mathcal{F}^Y_{1,n}$ is just

$$(65) \qquad \left(\frac{dQ}{dP}\right)_{\mathcal{F}^Y_{1,n}} = \frac{q(Y_1, Y_2, \ldots, Y_n)}{p(Y_1, Y_2, \ldots, Y_n)},$$

where

$$q(y_1, y_2, \ldots, y_n) := Q\{Y_j = y_j \ \forall \, 1 \leq j \leq n\}$$

and

$$p(y_1, y_2, \ldots, y_n) := P\{Y_j = y_j \ \forall \, 1 \leq j \leq n\}$$
$$= p^{\sum_{j=1}^n y_j} (1-p)^{n - \sum_{j=1}^n y_j}.$$

The Radon–Nikodym derivative on the larger $\sigma$-algebra $\mathcal{F}^{X,Y}_{1,n}$ cannot be so simply expressed, but is closely related to the Poissonized predictive probability $Z^*_{\alpha n}((\mathbf{X}, \mathbf{Y})_n)$ defined by (42). Define

$$(66) \qquad \hat{Z}_n := p(Y_1, Y_2, \ldots, Y_n)\left(\frac{dQ}{dP}\right)_{\mathcal{F}^{X,Y}_{1,n}};$$

then by (16) the random variable $\hat{Z}_{N(n)}$ has the same distribution under $P$ as does the Poissonized predictive probability (42) under $P_f$, for any $f$. Hence, the convergence (47) must also hold for the random variables $\hat{Z}_n$:

COROLLARY 20. *Under $P$, as $n \to \infty$,*

$$(67) \qquad n^{-1} \log \hat{Z}_n \xrightarrow{L^1} \psi_p(\alpha).$$

*Therefore, to prove the entropy inequality* (50), *it suffices to prove that*

$$(68) \qquad \lim_{n \to \infty} n^{-1} E_P \log \left( \frac{dQ}{dP} \right)_{\mathcal{F}_{1,n}^Y} < 0.$$

PROOF. The first assertion follows directly from (47) of Proposition 11. Thus, to prove the entropy inequality $\psi_p(\alpha) < -H(p)$, it suffices, in view of (66), to prove that (68) holds when the $\sigma$-algebra $\mathcal{F}_{1,n}^Y$ is replaced by $\mathcal{F}_{1,n}^{X,Y}$. But the former is a sub-$\sigma$-algebra of the latter; since log is a concave function, Jensen's inequality implies that

$$E_P \log \left( \frac{dQ}{dP} \right)_{\mathcal{F}_{1,n}^{X,Y}} \leq E_P \log \left( \frac{dQ}{dP} \right)_{\mathcal{F}_{1,n}^Y}. \qquad \square$$

5.2. *Digression*: *The relative SMB theorem.* The existence of the limit (67) is closely related to the *relative Shannon–McMillan–Breiman theorem* studied by several authors [14, 15, 18, 19]. The sequence $Y_1, Y_2, \ldots$ is, under either measure $P$ or $Q$, an ergodic stationary sequence of Bernoulli random variables. Thus, by the usual Shannon–MacMillan–Breiman theorem [27], as $n \to \infty$,

$$n^{-1} \log q(Y_1, Y_2, \ldots, Y_n) \xrightarrow{\text{a.s. } Q} -h_Q$$

and

$$n^{-1} \log p(Y_1, Y_2, \ldots, Y_n) \xrightarrow{\text{a.s. } P} -H(p),$$

where $h_Q$ is the Kolmogorov–Sinai entropy of the sequence $Y_j$ under $Q$. In general, of course, the almost sure convergence holds only for the probability measure indicated—see, for instance, [14] for an example where the first convergence fails under the alternative measure $P$. The *relative* Shannon–MacMillan–Breiman theorem of [19] gives conditions under which the difference of the two averages

$$n^{-1} \log \frac{q(Y_1, Y_2, \ldots, Y_n)}{p(Y_1, Y_2, \ldots, Y_n)} = n^{-1} \log \left( \frac{dQ}{dP} \right)_{\mathcal{F}_{1,n}^Y}$$

converges under $P$. In the case at hand, unfortunately, these conditions are not of much use: they essentially require the user to verify that

$$n^{-1} q(Y_1, Y_2, \ldots, Y_n) \xrightarrow{\text{a.s. } P} C$$

for some constant $C$. Thus, it appears that [19] does not provide a shortcut to the convergence (47).

5.3. *The rechargeable Pólya urn.* In the ordinary Pólya urn scheme, balls are drawn at random from an urn, one at a time; after each draw, the ball drawn is returned to the urn along with another of the same color. If initially the urn contains one red and one blue ball, then the limiting fraction $\Theta$ of red balls is uniformly distributed on the unit interval. The Pólya urn is connected with Bayesian statistics in the following way: the *conditional* distribution of the sequence of draws *given* the value of $\Theta$ is that of i.i.d. Bernoulli-$\Theta$ random variables.

The *rechargeable* Pólya urn is a simple variant of the scheme described above, differing only in that, before each draw, with probability $r > 0$, the urn is emptied and then reseeded with one red and one blue ball. Unlike the usual Pólya urn, the rechargeable Pólya urn is recurrent, that is, if $V_n := (R_n, B_n)$ denotes the composition of the urn after $n$ draws, then $V_n$ is a *positive recurrent* Markov chain on the state space $\mathbb{N} \times \mathbb{N}$. Consequently, $\{V_n\}$ may be extended to $n \in \mathbb{Z}$ in such a way that the resulting process is stationary. Let $Y_n$ denote the binary sequence recording the results of the successive draws ($1 = $ BLUE, $0 = $ RED). Clearly, this sequence has the same law as does the sequence $Y_1, Y_2, \ldots$ under the probability measure $Q$ [with $r = \alpha/(1 + \alpha)$].

LEMMA 21. *For any $\varepsilon > 0$, there exists $m$ such that the following is true: For any finite sequence $y_{-k}, y_{-k-1}, \ldots, y_0$, the conditional distribution of $Y_{m+1}, Y_{m+2}, \ldots, Y_{2m}$ given that $Y_i = y_i$ for all $-k \le i \le 0$ differs from the $Q$-unconditional distribution by less than $\varepsilon$ in total variation norm.*

PROOF. It is enough to show that the conditional distribution of $Y_{m+1}, \ldots, Y_{2m}$ given the composition $V_0$ of the urn before the first draw differs from the unconditional distribution by less than $\varepsilon$. Let $T$ be the time of the first regeneration (emptying of the urn) after time 0; then conditional on $T = n$, for any $n \le m$, and on $V_0$, the distribution of $Y_{m+1}, \ldots, Y_{2m}$ does not depend on the value of $V_0$. Thus, if $m$ is sufficiently large that the probability of having at least one regeneration event between the first and $m$th draws exceeds $1 - \varepsilon$, then the conditional distribution given $V_0$ differs from the unconditional distribution by less than $\varepsilon$ in total variation norm. $\square$

The construction of the sequence $\mathbf{Y} = Y_1, Y_2, \ldots$ using the rechargeable Pólya urn shows that this sequence behaves as a "factor" of a denumerable-state Markov chain (in terminology more familiar to statisticians, the sequence $\mathbf{Y}$ follows a "hidden Markov model"). Note that the original specification of the measure $Q$, in Section 5.1 above, exhibits $\mathbf{Y}$ as a factor of the Harris-recurrent Markov chain obtained by adjoining to the state variable the current value of $\mathbf{W}$. It does not appear that $Y_n$ can be represented as a function of a *finite*-state Markov chain; if it could, then results of Kaijser [12] would imply the existence of the limit

$$\lim_{n\to\infty} n^{-1}\log q(Y_1, Y_2, \ldots, Y_n)$$

almost surely under $P$, and exhibit it as the top Lyapunov exponent of a sequence of random matrix products. Unfortunately, little is known about the asymptotic behavior of random *operator* products (see [13] and references therein for the state of the art), and so it does not appear that (65) can be obtained by an infinite-state extension of Kaijser's result.

5.4. *Proof of* (68). Since it is not necessary to establish the convergence of the integrands on the left-hand side of (68), we shall not attempt to do so. Instead, we will proceed from the identity

$$(69)\quad n^{-1}E_P\log\frac{q(Y_1, Y_2, \ldots, Y_n)}{p(Y_1, Y_2, \ldots, Y_n)} = n^{-1}\sum_{k=0}^{n-1}E_P\log\frac{q(Y_{k+1}|Y_1, Y_2, \ldots, Y_k)}{p(Y_{k+1}|Y_1, Y_2, \ldots, Y_k)}.$$

Because the random variables $Y_i$ are i.i.d. Bernoulli-$p$ under $P$, the conditional probabilities $p(y_{k+1}|y_1, y_2, \ldots, y_k)$ must coincide with the *unconditional* probabilities $p(y_{k+1})$. Thus, the usual information inequality (Jensen's inequality), in the form $E_f\log(g(X)/f(X)) < 0$ for distinct probability densities $f, g$, implies that, for each $k$,

$$(70)\qquad E_P\log\frac{q(Y_{k+1}|Y_1, Y_2, \ldots, Y_k)}{p(Y_{k+1})} \le 0,$$

with the inequality strict unless the $Q$-conditional distribution of $Y_{k+1}$ given the past coincides with the Bernoulli-$p$ distribution. Moreover, the left-hand side of (70) will remain bounded away from 0 as long as the conditional distribution remains bounded away from the Bernoulli-$p$ distribution (in any reasonable metric, e.g., the total variation distance). Thus, to complete the proof of (68), it suffices to establish the following lemma.

LEMMA 22. *There is no sequence of integers $k_n \to \infty$ along which*

$$(71)\qquad \|q(\cdot|Y_1, Y_2, \ldots, Y_{k_n}) - p(\cdot)\|_{\mathrm{TV}} \longrightarrow 0$$

*in $P$-probability.*

PROOF. This is based on the fact that the sequence of draws $Y_1, Y_2, \ldots$ produced by the rechargeable Pólya urn is *not* a Bernoulli sequence, that is, the $Q$- and $P$-distributions of the sequence $Y_1, Y_2, \ldots$ are distinct. Denote by $q_k$ the $Q$-conditional probability that $Y_{k+1} = 1$ given the values $Y_1, Y_2, \ldots, Y_k$. Suppose that $q_{k_n} \to p$ in $P$-probability; then by summing over successive values of the last $l$ variables, it follows that $q_{k_n-l} \to p$ in $P$-probability for each fixed $l \in \mathbb{N}$. We will show that this leads to a contradiction.

Consider the following method of generating binary random variables $Y_1, Y_2, \ldots, Y_{2m}$: first generate i.i.d. Bernoulli-$p$ random variables $Y_j$ for $-k \leq j \leq 0$; then, conditional on their values, generate $Y_1$ according to $q_{k+1}$; then, conditional on $Y_1$, generate $Y_2$ according to $q_{k+2}$; and so on. By the hypothesis of the preceding paragraph, there is a sequence $k_n \to \infty$ such that, for any fixed $m$, the joint distribution of $Y_1, Y_2, \ldots, Y_{2m}$ converges to the product-Bernoulli-$p$ distribution. But this contradicts the mixing property of the rechargeable Pólya urn asserted by Lemma 21 above. $\quad\square$

## APPENDIX A: AN ALMOST SUBADDITIVE WLLN

The purpose of this appendix is to prove the simple variant of the subadditive ergodic theorem required in Section 4. For the original subadditive ergodic theorem of Kingman, see [16], and for another variant that is useful in applications to percolation theory, see [17]. There are two novelties in our version: (a) the subadditivity relation is only approximate, with a random error; and (b) there is no measure-preserving transformation related to the sequence $S_n$.

PROPOSITION A.1.   *Let $S_n$ be real random variables. Suppose that, for each pair $m, n \geq 1$ of positive integers, there exist random variables $S'_{m,m+n}$, $S''_{n,m+n}$ and a nonnegative random variable $R_{m,n}$ such that:*

(a) *$S'_{m,m+n}$ and $S''_{n,m+n}$ are independent;*
(b) *$S'_{m,m+n}$ has the same distribution as $S_m$;*
(c) *$S''_{m,m+n}$ has the same distribution as $S_n$;*
(d) *the random variables $\{R_{m,n}\}_{m,n \geq 1}$ are identically distributed;*
(e) *$ER_{1,1} < \infty$ and $\{S_n/n\}_{n \geq 1}$ are uniformly integrable; and*
(f) *for all $m, n \geq 1$,*

(A.1) $$S_{m+n} \leq S'_{m,m+n} + S''_{n,m+n} + R_{m,n}.$$

*Then*

(A.2) $$\frac{S_n}{n} \xrightarrow{L^1} \gamma := \liminf_{n \to \infty} \frac{ES_n}{n}.$$

NOTE.   The random variables $\{S_n/n\}$ considered in Corollary 15 are uniformly bounded, and so the uniform integrability hypothesis (e) holds trivially.

PROOF OF PROPOSITION A.1.   Since the random variables $S'_{m,m+n}$ and $S''_{n,m+n}$ are independent, with the same distributions as $S_m$ and $S_n$, respectively, Carathéodory's theorem on extension of measures implies that the probability space may be enlarged so as to support additional random variables permitting recursion on the inequality (A.1). Here the simplest recursive strategy works: from a starting value $n = km + r$, reduce by $m$ at each step. This leads to an inequality of the form

$$(A.3) \qquad S_{km+r} \leq S_r^0 + \sum_{j=1}^k S_m^j + \sum_{j=1}^k R_j,$$

where the random variables $\{S_m^j\}_{j \geq 1}$ are i.i.d., each with the same distribution as $S_m$, and the random variables $R_j$ are identically distributed (but not necessarily independent), each with the law of $R_{1,1} := R$.

The weak law (A.2) is easily deduced from the inequality (A.3). Note first that the special case of (A.1) with $m = 1$, together with hypothesis (d), implies that $ES_n \leq nER + nES_1 < \infty$ for every $n \geq 1$, and so $\gamma < \infty$. Assume for definiteness that $\gamma > -\infty$; the case $\gamma = -\infty$ may be treated by a similar argument. Divide each side of (A.3) by $km$; as $k \to \infty$,

$$\frac{S_r^0}{km} \xrightarrow{P} 0 \quad \text{and} \quad \frac{1}{km} \sum_{j=1}^k S_m^j \xrightarrow{P} \frac{ES_m}{m},$$

the latter by the usual WLLN. The WLLN need not apply to the sum $\sum R_j$, since the terms are not necessarily independent; however, since all of the terms are nonnegative and have the same expectation $ER < \infty$, Markov's inequality implies that, for any $\varepsilon > 0$,

$$P\left\{\frac{1}{km} \sum_{j=1}^k R_j \geq \varepsilon\right\} \leq \frac{ER}{m\varepsilon}.$$

Thus, letting $m \to \infty$ through a subsequence along which $ES_m/m \to \gamma$, we find that, for any $\varepsilon > 0$,

$$\lim_{n \to \infty} P\{S_n \geq n\gamma + n\varepsilon\} = 0.$$

Since the r.v.s $S_n/n$ are uniformly integrable, this implies that

$$\lim_{n \to \infty} P\{S_n \leq n\gamma - n\varepsilon\} = 0,$$

because otherwise $\liminf ES_n/n < -\gamma$. This proves that $S_n/n \to \gamma$ in probability; in view of the uniform integrability of the sequence $S_n/n$, convergence in $L^1$ follows. $\square$

REMARK. Numerous variants of this proposition are true, and may be established by more careful recursions. Among these are SLLNs for random variables satisfying hypotheses such as those given in (59) above, where the remainders $\log \Lambda(m + n)$ are not identically distributed, but whose growth is sublinear in $m + n$. For hints as to how such results may be approached, see [11].

## APPENDIX B: PROOF OF PROPOSITION 2

Recall that an interval $J$ is $\varepsilon$-bad if any one of the inequalities (30), (31) or (32) holds. For each of these inequalities, there are two possibilities: the relevant count $N(J), N^S(J), N^F(J)$ may be unusually large or unusually small. Thus, there are six distinct ways that $J$ may be $\varepsilon$-bad, and hence six ways that a point $x$ may be $(\varepsilon, \kappa)$-bad. To prove (33), we will partition the set $\mathcal{B} = \mathcal{B}_n(\varepsilon, \kappa)$ into six subsets, one for each possibility, and show that (33) holds for each of the six subsets. In fact, since the six inequalities (30)–(32) are coupled [in particular, if $N^S(J)$ is unusually small, then either $N(J)$ is also unusually small or $N^F(J)$ is unusually large], it suffices to consider only four possibilities: those where $N(J)$ is either unusually large or small, and those where $N^S(J)$ or $N^F(J)$ is unusually large. These may all be handled in a similar fashion, so we shall consider only the possibilities involving large discrepancies of $N(J)$. Thus, set

$$\mathcal{B}^+ = \left\{ x \colon \sup_{J \colon x \in J; |J| \geq \kappa/n} N(J)/|J| \geq (1+\varepsilon)n \right\},$$

$$\mathcal{B}^- = \left\{ x \colon \sup_{J \colon x \in J; |J| \geq \kappa/n} N(J)/|J| \leq (1-\varepsilon)n \right\}.$$

We will show that, for any $\varepsilon > 0$, there are positive constants $\kappa, \gamma, C$ such that

$$(\text{B.1}) \qquad P_f\{|\mathcal{B}^\pm| \geq \varepsilon\} \leq Ce^{-\gamma n}.$$

The proof of (B.1) is of a familiar type in the theory of empirical processes: see, for instance, [20], Chapter 3 for related arguments. The strategy is to *bracket* each of the bad sets $\mathcal{B}^\pm$ by nearby sets in finite $\sigma$-algebras whose cardinalities are small compared to $e^n$. To carry out this bracketing, we will call on a weak form of the Vitali covering lemma (see [24], Chapter 1, Lemma 1.6):

COVERING LEMMA.   *Let $F$ be a measurable subset of $[0,1]$ that is covered by a collection $\mathcal{V}$ of subintervals of $[0,1]$. Then there exist pairwise disjoint intervals $I_1, I_2, \ldots$ in $\mathcal{V}$ such that*

$$(\text{B.2}) \qquad \sum_j |I_j| \geq \tfrac{1}{5}|F|.$$

For definiteness, consider the set $\mathcal{B}^+$. This set is, by definition, the union of intervals $J$ of lengths $\geq \kappa/n$ (not necessarily pairwise disjoint!), each satisfying

$$(\text{B.3}) \qquad N(J) \geq (1+\varepsilon)n|J|.$$

By the Covering lemma, there is a collection $\mathcal{J}$ of pairwise disjoint intervals $J_i$ among these whose lengths sum to at least $|\mathcal{B}^+|/5$. Because the lengths of these intervals are bounded below by $\kappa/n$, there are at most $n/\kappa$ intervals in $\mathcal{J}$.

Let $m = m_n$ be the smallest integer such that $m > 8n/(\varepsilon\kappa)$. For each interval $J_i \in \mathcal{J}$, let $J_i'$ be the minimal interval containing $J_i$ with endpoints of the form $j/m$, where $j \in \mathbb{Z}$. Observe that the intervals $J_i'$ need not be pairwise disjoint, but keep in mind that the intervals $J_i$ are. Moreover, $|J_i| \leq |J_i'| \leq |J_i| + 2/m$, so by (B.3),

$$\begin{aligned}
N(J_i) &\geq (1+\varepsilon)n(|J_i'| - 2/m) \\
&\geq (1+\varepsilon')n|J_i'|,
\end{aligned} \tag{B.4}$$

where $\varepsilon' = \varepsilon/2$. Define $\mathcal{B}^*$ to be the union of the intervals $J_i'$. By construction, $\mathcal{B}^*$ is a union of intervals $[j/m, (j+1)/m]$ and contains $\bigcup_{\mathcal{J}} J_i$; thus, by (B.4), since the intervals $J_i \in \mathcal{J}$ are pairwise disjoint,

$$N(\mathcal{B}^*) \geq \sum_{\mathcal{J}} N(J_i) \geq (1+\varepsilon')n|\mathcal{B}^*|. \tag{B.5}$$

Now recall that the collection $\mathcal{J}$ was chosen so that $\sum_{\mathcal{J}} |J_i| \geq |\mathcal{B}^+|/5$. Since $\mathcal{B}^*$ contains $\bigcup_{\mathcal{J}} J_i$, it follows that, on the event $\{|\mathcal{B}^+| \geq \varepsilon\}$,

$$|\mathcal{B}^*| \geq \varepsilon/5. \tag{B.6}$$

Hence, to bound the probability (B.1), it suffices to bound the probability that inequality (B.6) obtains. Observe that there are precisely $2^m$ possibilities for the set $\mathcal{B}^*$. Let $B$ be such a possibility, and suppose that $|B| \geq \varepsilon/5$. Under $P_f$, the count $N(B)$ has the Binomial distribution with parameters $n, |B|$. By a standard concentration inequality for the Binomial distribution (e.g., Hoeffding's inequality), there exist constants $C > 0$ and $\rho = \rho(\varepsilon, \varepsilon') > 0$ such that

$$P_f\{N(B) \geq (1+\varepsilon')n|B|\} \leq Ce^{-\rho n}.$$

Therefore,

$$P_f\{|\mathcal{B}^*| \geq \varepsilon\} \leq C2^m e^{-\rho n}.$$

Finally, recall that $m = \lceil 8n/(\varepsilon\kappa) \rceil$. Thus, if $\varepsilon\kappa$ is sufficiently large, then

$$2^m e^{-\rho n} \leq e^{-\gamma n}$$

for some $\gamma > 0$, and so (B.1) follows for $\mathcal{B}^+$. A similar argument (with bracketing from the inside rather than from the outside) applies for $\mathcal{B}^-$.

**Acknowledgment.** We thank Peter Radchenko for carefully reading the first version of this article and suggesting a number of improvements in the exposition.

## REFERENCES

[1] BARRON, A., SCHERVISH, M. J. and WASSERMAN, L. (1999). The consistency of posterior distributions in nonparametric problems. *Ann. Statist.* **27** 536–561. MR1714718

[2] CHI, Z. (2001). Stochastic sub-additivity approach to the conditional large deviation principle. *Ann. Probab.* **29** 1303–1328. MR1872744

[3] CORAM, M. (2002). Nonparametric Bayesian classification. Ph.D. dissertation, Stanford Univ.

[4] DIACONIS, P. and FREEDMAN, D. (1986). On inconsistent Bayes estimates of location. *Ann. Statist.* **14** 68–87. MR0829556

[5] DIACONIS, P. and FREEDMAN, D. A. (1993). Nonparametric binary regression: A Bayesian approach. *Ann. Statist.* **21** 2108–2137. MR1245784

[6] DIACONIS, P. and FREEDMAN, D. A. (1995). Nonparametric binary regression with random covariates. *Probab. Math. Statist.* **15** 243–273. MR1369802

[7] FREEDMAN, D. A. (1963). On the asymptotic behavior of Bayes' estimates in the discrete case. *Ann. Math. Statist.* **34** 1386–1403. MR0158483

[8] FREEDMAN, D. and DIACONIS, P. (1983). On inconsistent Bayes estimates in the discrete case. *Ann. Statist.* **11** 1109–1118. MR0720257

[9] GHOSAL, S., GHOSH, J. K. and VAN DER VAART, A. W. (2000). Convergence rates of posterior distributions. *Ann. Statist.* **28** 500–531. MR1790007

[10] GHOSH, J. K. and RAMAMOORTHI, R. V. (2003). *Bayesian Nonparametrics.* Springer, New York. MR1992245

[11] HAMMERSLEY, J. M. (1962). Generalization of the fundamental theorem on sub-additive functions. *Proc. Cambridge Philos. Soc.* **58** 235–238. MR0137800

[12] KAIJSER, T. (1975). A limit theorem for partially observed Markov chains. *Ann. Probab.* **3** 677–696. MR0383536

[13] KARLSSON, A. and MARGULIS, G. A. (1999). A multiplicative ergodic theorem and nonpositively curved spaces. *Comm. Math. Phys.* **208** 107–123. MR1729880

[14] KIEFFER, J. C. (1973). A counterexample to Perez's generalization of the Shannon–McMillan theorem. *Ann. Probab.* **1** 362–364. MR0351626

[15] KIEFFER, J. C. (1974). A simple proof of the Moy–Perez generalization of the Shannon–McMillan theorem. *Pacific J. Math.* **51** 203–206. MR0347448

[16] KINGMAN, J. F. C. (1973). Subadditive ergodic theory (with discussion). *Ann. Probab.* **1** 883–909. MR0356192

[17] LIGGETT, T. M. (1985). An improved subadditive ergodic theorem. *Ann. Probab.* **13** 1279–1285. MR0806224

[18] PEREZ, A. (1964). Extensions of Shannon–McMillan's limit theorem to more general stochastic processes. In *Trans. Third Prague Conference on Information Theory, Statistical Decision Functions, Random Processes* (*Liblice, 1962*) 545–574. Publ. House Czech. Acad. Sci., Prague. MR0165996

[19] PEREZ, A. (1980). On Shannon–McMillan's limit theorem for pairs of stationary random processes. *Kybernetika* (*Prague*) **16** 301–314. MR0591959

[20] POLLARD, D. (1984). *Convergence of Stochastic Processes.* Springer, New York. MR0762984

[21] RUELLE, D. (1999). *Statistical Mechanics. Rigorous Results.* World Scientific, River Edge, NJ. MR1747792

[22] SCHWARTZ, L. (1965). On Bayes procedures. *Z. Wahrsch. Verw. Gebiete* **4** 10–26. MR0184378

[23] SHEN, X. and WASSERMAN, L. (2001). Rates of convergence of posterior distributions. *Ann. Statist.* **29** 687–714. MR1865337

[24] STEIN, E. M. (1970). *Singular Integrals and Differentiability Properties of Functions.* Princeton Univ. Press. MR0290095

[25] WALKER, S. G. (2004). New approaches to Bayesian consistency. *Ann. Statist.* **32** 2028–2043. MR2102501

[26] WALKER, S. G. (2004). Modern Bayesian asymptotics. *Statist. Sci.* **19** 111–117. MR2082150

[27] WALTERS, P. (1982). *An Introduction to Ergodic Theory.* Springer, New York. MR0648108

DEPARTMENT OF STATISTICS
UNIVERSITY OF CHICAGO
5734 UNIVERSITY AVENUE
CHICAGO, ILLINOIS 60637
USA
E-MAIL: coram@galton.uchicago.edu
        lalley@galton.uchicago.edu